\documentclass[12pt]{report}

\usepackage{fancyhdr,mathptmx}

\fancypagestyle{plain}

\usepackage{pictexwd}
\usepackage{amsmath}
\usepackage{amscd}
\usepackage{amssymb}
\usepackage{amsthm}
\usepackage{xspace}
\usepackage[mathscr]{eucal}
\usepackage[all,tips]{xy}
\usepackage[dvips]{graphicx}
\usepackage{fancyhdr}
\usepackage{verbatim}
\usepackage{syntonly}

\theoremstyle{plain}

\newtheorem{thm}{Theorem}[section]
\newtheorem*{thmnl}{Theorem}
\newtheorem{cor}[thm]{Corollary}

\newtheorem{prop}[thm]{Proposition}

\newtheorem{lemma}[thm]{Lemma}

\theoremstyle{definition}

\newtheorem*{sol}{Solution}

\newtheorem{exa}{Example}[section]

\newtheorem*{prf}{Proof}
\newtheorem*{rem}{Remark}

\newenvironment{pf}
{\begin{proof}} {\end{proof}}

\DeclareMathOperator{\Aut}{Aut} 
\DeclareMathOperator{\Inn}{Inn}
\DeclareMathOperator{\Out}{Out}
\DeclareMathOperator{\Stab}{Stab}

\DeclareMathOperator{\rank}{rank}
\DeclareMathOperator{\id}{id}

\DeclareMathOperator*{\Ima}{Im} 
\DeclareMathOperator*{\RE}{Re}

\DeclareMathOperator*{\ch}{char} 
\DeclareMathOperator{\tr}{Tr}
\DeclareMathOperator{\Gal}{Gal}

\DeclareMathOperator{\SL}{SL}
\DeclareMathOperator{\PSL}{PSL}
\DeclareMathOperator{\GL}{GL}

\DeclareMathOperator{\Ort}{O}
\DeclareMathOperator{\PO}{PO}

\DeclareMathOperator{\SO}{SO}

\DeclareMathOperator{\Uni}{U}
\DeclareMathOperator{\PU}{PU}

\DeclareMathOperator{\SU}{SU}

\DeclareMathOperator{\Sp}{Sp}

\DeclareMathOperator{\PSp}{PSp}

\DeclareMathOperator{\Sol}{Sol}
\DeclareMathOperator{\Nil}{Nil}

\DeclareMathOperator{\Mat}{M}
\DeclareMathOperator{\Aff}{Aff}
\DeclareMathOperator{\Sim}{Sim}
\DeclareMathOperator{\Euc}{Euc}

\DeclareMathOperator*{\vol}{vol}

\DeclareMathOperator*{\Isom}{Isom}

\DeclareMathOperator*{\sign}{sign}

\newcommand{\eps}{\varepsilon}
\newcommand{\vp}{\varphi}

\newcommand{\al}{\alpha}

\newcommand{\be}{\beta}

\newcommand{\ga}{\gamma}
\newcommand{\Ga}{\Gamma}

\newcommand{\te}{\theta}

\newcommand{\si}{\sigma}

\newcommand{\om}{\omega}

\newcommand{\de}{\delta}

\newcommand{\la}{\lambda}
\newcommand{\La}{\Lambda}

\newcommand{\ol}{\overline}
\newcommand{\nl}{\ensuremath{\newline}}

\newcommand{\ig}{\includegraphics}

\newcommand{\wt}{\widetilde}
\newcommand{\wh}{\widehat}

\newcommand{\iny}{\infty}

\newcommand{\tri}{\ensuremath{\triangle}}
\newcommand{\prt}{\partial}

\newcommand{\co}{\ensuremath{\colon}}

\newcommand{\innp}[1]{\left< #1 \right>}
\newcommand{\abs}[1]{\left\vert#1\right\vert}
\newcommand{\set}[1]{\left\{#1\right\}}
\newcommand{\brac}[1]{\left[#1\right]}
\newcommand{\pr}[1]{\left( #1 \right) }

\newcommand{\norm}[1]{\left\Vert#1\right\Vert}

\newcommand{\su}{\subset}

\newcommand{\bu}{\bigcup}

\newcommand{\op}{\oplus}

\newcommand{\bop}{\bigoplus}

\newcommand{\smin}{\setminus}

\newcommand{\bdef}{\overset{\text{def}}{=}}

\newcommand{\lra}{\longrightarrow}

\newcommand{\lmto}{\longmapsto}

\newcommand{\B}[1]{\ensuremath{\mathbf{#1}}}
\newcommand{\BB}[1]{\ensuremath{\mathbb{#1}}}
\newcommand{\Cal}[1]{\ensuremath{\mathcal{#1}}}
\newcommand{\Fr}[1]{\ensuremath{\mathfrak{#1}}}

\newcommand{\Hy}{\ensuremath{\B{H}}}
\newcommand{\N}{\ensuremath{\B{N}}}
\newcommand{\Q}{\ensuremath{\B{Q}}}
\newcommand{\R}{\ensuremath{\B{R}}}
\newcommand{\Z}{\ensuremath{\B{Z}}}
\newcommand{\C}{\ensuremath{\B{C}}}

\newcommand{\res}{\ensuremath{\textrm{Res}}}

\newcommand{\refP}[1]{Proposition~\ref{P:#1}}

\newcommand{\refT}[1]{Theorem~\ref{T:#1}}
\newcommand{\refL}[1]{Lemma~\ref{L:#1}}

\newcommand{\refC}[1]{Corollary~\ref{C:#1}}
\newcommand{\refE}[1]{Equation~(\ref{E:#1})}

\begin{document}

\author{David Ben McReynolds}
\title{Cusps of arithmetic orbifolds}
\maketitle

\begin{abstract}
This thesis investigates cusp cross-sections of arithmetic real, complex, and quaternionic hyperbolic $n$--orbifolds. We give a smooth classification of these submanifolds and analyze their induced geometry. One of the primary tools is a new subgroup separability result for general arithmetic lattices.
\end{abstract}

\tableofcontents 
\listoffigures

\chapter{Introduction and main results}

\section{Bounding and geometric bounding questions}

It is a classical problem in topology to decide whether or not a closed $n$--manifold $M$ bounds. Hamrick and Royster \cite{HamrickRoyster82} resolved this in the affirmative for flat $n$--manifolds and Rohlin \cite{Rohlin51} for closed 3--manifolds (see also Rourke \cite{Rourke85} and Milnor \cite{Milnor63}). However, beyond these two classes there are few other settings where the story is nearly this complete.

The introduction of geometry to a topological problem provides additional structure which can lead to new insight. This philosophy serves as motivation for the primary concern of this thesis which is a geometric notion of bounding and its specialization to flat and almost flat manifolds. To both guide the reader and establish some basic terminology required for our results, we review the simplest case of this venture whose concern is with the class of finite volume hyperbolic $n$--manifolds and their flat cusp cross-sections.

Every finite volume hyperbolic $n$--orbifold $W$ has a thick-thin decomposition comprised of a compact manifold $W_{core}$ with boundary components $F_1,\dots,F_m$ and manifolds $E_1,\dots,E_m$ of the form $F_j \times \R^+$. The manifolds $E_j$ are called \emph{cusp ends}, the manifolds $F_j$ are called \emph{cusp cross-sections}, and the union of $W_{core}$ with $E_1,\dots,E_m$ recovers $W$ topologically. The inclusion $F_j \lra W$ induces an injective homomorphism $\pi_1(F_j) \lra \Isom(\Hy^n)$ whose image is virtually contained in a subgroup isomorphic to $\R^{n-1}$. Consequently, $\pi_1(F_j)$ is virtually abelian, $F_j$ is a flat orbifold of dimension $n-1$, and each $F_j$ is a totally geodesic boundary submanifold of $W_{core}$.

\begin{figure}[h!]
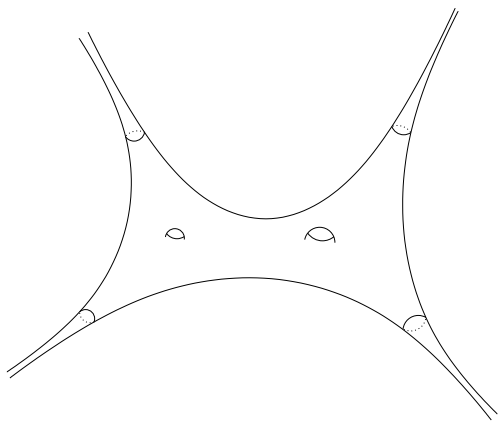

\begin{center}
\ig{CuspedManifold}
\caption{A finite volume hyperbolic surface with four cusp ends.}
\end{center}
\end{figure}

We say a flat $n$--manifold $F$ \emph{arises as a cusp cross-section} if there exists a hyperbolic $(n+1)$--manifold $W$ and a boundary component $F_j$ of $W_{core}$ diffeomorphic to $F$. More generally, we say $M$ \emph{geometrically bounds} if there exists a compact Riemannian $(n+1)$--manifold $W$ with totally geodesic boundary $\prt W$ and a diffeomorphism between $\prt W$ and $M$.

Farrell and Zdravkovska \cite{FarrellZdravkovska83} conjectured every flat $n$--manifold arises as a cusp cross-section of a 1--cusped hyperbolic $(n+1)$--manifold and this is easily verified for $n=2$. Indeed, the complement of a knot in $S^3$ is typically endowed with a finite volume, complete hyperbolic structure with one cusp (see \cite{Thurston82}), and thus gives the realization of the 2--torus $T^2$ as a cusp cross-section of a 1--cusped hyperbolic 3--manifold. Likewise, the Klein bottle arises as a cusp cross-section of the 1--cusped Gieseking manifold (see \cite{Ratcliffe94}). However, Long and Reid \cite{LongReid00} constructed counterexamples by showing any flat 3--manifold arising as a cusp cross-section of a 1--cusped hyperbolic 4--manifold must have integral $\eta$--invariant.\smallskip\smallskip

The failure of the conjecture of Farrell--Zdravkovska is far from total. Nimershiem \cite{Nimershiem98} showed every flat 3--manifold arises as a cusp cross-section of a hyperbolic 4--manifold, and Long and Reid \cite{LongReid02} proved every flat $n$--manifold arises as a cusp cross-section of a hyperbolic $(n+1)$--orbifold. 

Our main aim is a complete understanding of the topological and geometric structure of cusp cross-sections of hyperbolic $n$--orbifolds and their complex and quaternionic analogs. The absence of a general geometric construction for real, complex, and quaternionic hyperbolic orbifolds forces our restriction to orbifolds produced by arithmetic means. Given this forced restriction, the picture we provide here is nearly complete.  

\section{Classifying cusp cross sections of $X$--hyperbolic orbifolds}

For $X = \R$, $\C$, or $\BB{H}$, cusp cross-sections of finite volume $X$--hyperbolic $(n+1)$--orbifolds are flat $n$--manifolds or almost flat orbifolds modelled on the $(2n+1)$--dimensional Heisenberg group $\Fr{N}_{2n+1}$ or the $(4n+3)$--dimensional quaternionic Heisenberg group $\Fr{N}_{4n+3}(\BB{H})$. Our first result shows the result of Long--Reid \cite{LongReid02} does not extend to the complex or quaternionic hyperbolic settings. \footnote{Since every lattice in the isometry group of quaternionic hyperbolic space is arithmetic (see \cite{Corlette92} and \cite{GromovSchoen92}), an arithmeticity assumption in (b) is superfluous.} 

\begin{thm}\label{T:Ob}
\begin{description}
\item[(a)] For every $n \geq 2$, there exist infinite families of closed almost flat $(2n+1)$--manifolds modelled on $\Fr{N}_{2n+1}$ that are not diffeomorphic to a cusp cross-section of any arithmetic complex hyperbolic $(n+1)$--orbifold. 
\item[(b)] For every $n \geq 1$, there exist infinite families of closed almost flat $(4n+3)$--manifolds modelled on $\Fr{N}_{4n+3}(\BB{H})$ that are not diffeomorphic to a cusp cross-section of any finite volume quaternionic hyperbolic $(n+1)$--orbifold. 
\end{description}
\end{thm}

\refT{Ob} is a consequence of the topological classification of cusp cross-sections of arithmetic $X$--hyperbolic $(n+1)$--orbifolds (see \refT{324}). Despite the negative nature of \refT{Ob}, it is misleading. There are a wide range of manifolds which do arise and the metrics present vary considerably. One such example is the special case of complex hyperbolic 2--orbifolds and their $\Nil$ 3--manifold cusp cross-sections.

\begin{thm}\label{T:Nil3Thm}
Every $\Nil$ 3--manifold is diffeomorphic to a cusp cross-section of an arithmetic complex hyperbolic 2--orbifold.
\end{thm}

\section{Peripheral subgroup separability}

In order to give a complete classification of cusp cross-sections of arithmetic $X$--hyperbolic orbifolds, we require certain subgroup separability results. For a group $G$, a subgroup $H$ of $G$, and $g$ in $G \smin H$, we say $H$ and $g$ are \emph{separated} if there exists a subgroup $K$ of finite index in $G$ which contains $H$ but not $g$. We say $H$ is \emph{separable} in $G$ if for every $g$ in $G \smin H$, $g$ and $H$ are separated. We say $G$ is \emph{subgroup separable} (or LERF) if every finitely generated subgroup is separable.

The next theorem and its subsequent corollary provide one of the main tools in the classification of cusp cross-sections of arithmetic $X$--hyperbolic orbifolds.

\begin{thm}[Borel subgroup separability theorem]\label{T:Borel}
Let $k$ be a number field, $\B{G}$ a $k$--algebraic group, and $\B{B}$ a Borel subgroup of $\B{G}$. Then any subgroup of $\B{B}(\Cal{O}_k)$ is separable in $\B{G}(\Cal{O}_k)$. 
\end{thm}  

\begin{cor}[Stabilizer subgroup separability theorem]\label{C:Sep}
Let $Y$ be $\Hy_\R^n$, $\Hy_\C^n$, $\Hy_\BB{H}^n$ or $\Hy_\BB{O}^2$, $\La$ an arithmetic lattice in $\Isom(Y)$, and $v$ in $\prt Y$. Then every subgroup of $\La \cap \Stab(v)$ is separable in $\La$.    
\end{cor}

How these results are brought to bear on the topological problem of classifying cusp cross-sections of arithmetic $X$--hyperbolic $n$--orbifolds is fairly well known. This is achieved through a relationship established by Scott \cite{Scott78} between subgroup separability and the lifting in a finite cover of a $\pi_1$--injective immersion to a $\pi_1$--injective embedding. In tandem with \refC{Sep}, we obtain the following geometrically applicable corollary.

\begin{cor}\label{C:Immersion}
Let $\rho\co N \lra M$ be a $\pi_1$--injective immersion of an almost flat manifold $N$ into an arithmetic $Y$--hyperbolic $n$--orbifold. Then there exists a finite cover $\psi\co M^\prime \lra M$ such that $\rho$ lifts to an embedding.
\end{cor} 

\section{Density of cusp shapes}

It is well known that most similarity classes of flat structures on a flat $n$--manifold cannot arise in the cusp cross-sections of any hyperbolic $(n+1)$--orbifold. Nevertheless, Nimershiem \cite{Nimershiem98} showed for any flat 3--manifold, the similarity classes of flat metrics arising in the cusp cross-sections of hyperbolic 4--manifolds are dense in the space of flat similarity classes. She conjectured \cite[Conj. 2']{Nimershiem98} this for every flat $n$--manifolds, and our next result verifies this in the orbifold setting.

\begin{thm}\label{T:SDensity} For a flat $n$--manifold $M$, the flat similarity classes that arise in the cusp cross-sections of arithmetic hyperbolic $(n+1)$--orbifolds are dense in the space of flat similarity classes.
\end{thm}

Using Selberg's lemma, we verify the full conjecture for the $n$--torus.

\begin{cor}\label{C:Density}
For the $n$--torus, the flat similarity classes that arise in the cusp cross-sections of arithmetic hyperbolic $(n+1)$--manifolds are dense in the space of flat similarity classes.
\end{cor}

The similarity structures used in the proof of \refT{SDensity} are precisely those that arise in cusp cross-sections of arithmetic hyperbolic orbifolds. With an extension of this to infranil manifolds modelled on $\Fr{N}_{2n-1}$ and $\Fr{N}_{4n-1}(\BB{H})$, we obtain a geometric classification of cusp cross-sections of arithmetic $X$--hyperbolic $n$--orbifolds. 

\section{Hilbert modular varieties}

The same methods used for cusp cross-sections of $X$--hyperbolic orbifolds can be implemented for cusp cross-sections of Hilbert modular varieties over a totally real number field $k$. The corresponding cusps are virtual $n$--torus bundles over $(n-1)$--tori where $[k:\Q]=n$, $\rank \Cal{O}_k^\times = n-1$, and for brevity we call these \emph{virtual $(n,n-1)$--torus bundles}. We say $\be$ in $k$ is \emph{totally positive} if $\si_j(\be) >0$ for every real embedding $\si_j$ of $k$. We denote the set of totally positive elements and totally positive integers by $k_+$ and $\Cal{O}_{k,+}$, and define the sets $k_+^\times = k_+ \cap k^\times$, $\Cal{O}_{k,+}^\times = \Cal{O}_k^\times \cap \Cal{O}_{k,+}$. We say a virtual torus bundle $N$ is \emph{$k$--defined} if there exists a faithful representation
\[ \rho\colon \pi_1(N) \lra k \rtimes k_+^\times, \]
and we say $N$ is \emph{$k$--arithmetic} if $\rho(\pi_1(N))$ and $\Cal{O}_k \rtimes \Cal{O}_{k,+}^\times$ are commensurable.

\begin{thm}\label{T:Torus}
A virtual $(n,n-1)$--torus bundle $N$ is diffeomorphic to a cusp cross-section of a Hilbert modular variety over $k$ if and only if $\pi_1(N)$ is $k$--arithmetic.
\end{thm} 

\refT{Torus} answers a question of Hirzebruch \cite[page 203]{Hirzebruch73} who asked in our terminology which $k$--arithmetic torus bundles arise as cusp cross-sections of Hilbert modular varieties over $k$. \refT{Torus} also provides an obstruction to the realizability of a torus bundle as a cusp cross-section of a Hilbert modular variety, yielding $(n,n-1)$--torus bundles (with $n>2$) that are not diffeomorphic to a cusp cross-section of any (generalized) Hilbert modular variety.

Hilbert modular varieties over real quadratic number fields are traditionally called \emph{Hilbert modular surfaces} and possess $\Sol$ 3--manifold cusp cross-sections. Here, as with $\Nil$ 3--manifolds, we achieve the converse.

\begin{thm}\label{T:Sol3Thm}
Every $\Sol$ 3--manifold is diffeomorphic to a cusp cross-section of a generalized Hilbert modular surface.
\end{thm}

Our final result follows from an obstruction to geometric bounding whose derivation employs the methodology used by Long--Reid \cite{LongReid00} for flat 3--manifolds.

\begin{thm}\label{T:GB}
There exists a $\Sol$ 3--manifold that cannot be diffeomorphic to a cusp cross-section of any 1--cusped Hilbert modular surface with torsion free fundamental group.
\end{thm} 

\chapter{Preliminaries}

\section{Some basic algebra}

For a number field $k/\Q$, we denote the \emph{ring of $k$--integers} by $\Cal{O}_k$, the set of distinct real embeddings by $\si_1,\dots,\si_{r_1}$, and the set of distinct complex embeddings by $\tau_1,\dots,\tau_{r_2}$. We say that $k$ is \emph{totally real} if $r_2=0$, \emph{totally imaginary} if $r_1=0$, and call quadratic extensions $E/F$ with $E$ totally imaginary and $F$ totally real, \emph{CM fields}. 

Associated to each $\al \in k^\times$ and cyclic extension $L/k$ of degree $d$ with Galois group $\innp{\te}$,
is a central simple $k$--algebra $(L/k,\te,\al)$
\[ (L/k,\te,\al) = \set{\sum_{j=0}^{d-1} \be_j X^j~:~\be_j \in L,~X^d = \al,~X\be = \te(\be)X \text{ for }\be \in L} \]
called a \emph{cyclic algebra} of degree $d$ over $k$. In the event $L/k$ is quadratic, the resulting algebras are \emph{$k$--quaternion algebras} with \emph{Hamilton's quaternions} $\BB{H}$ being the most well known example. 

By a theorem of Wedderburn, any quaternion algebra over $\C$ is isomorphic to $\Mat(2;\C)$ and over $\R$ to either $\Mat(2;\R)$ or $\BB{H}$. For a $k$--quaternion algebra $B$, each embedding $\la$ of $k$ yields a new quaternion algebra $^\la B$. If  either $\la$ is complex or 
$^\la B \otimes_{\la(k)} \R \cong \Mat(2;\R)$, we say $^\la B$ is \emph{unramified} at $\la$, and otherwise say $^\la B$ is \emph{ramified} at $\la$. When $B$ is (un)ramified at each embedding $\la$ of $k$, we say $B$ is \emph{totally (un)ramified}.

By an \emph{$\Cal{O}_k$--order} in an $k$--algebra $A$, we mean a finitely generated subring $\Cal{O}$ of $A$ that is finitely generated as an $\Cal{O}_k$--module and $A = \Cal{O} \otimes_{\Cal{O}_k} k$. In the sequel we require an existence theorem for orders in central, simple $k$--algebras and refer the reader to \cite{Reiner03} for a proof.

\section{Algebraic groups and lattices}

For an algebraically closed field $\ol{k}$ of characteristic zero, a \emph{linear algebraic group} $\B{G}$ is a Zariski closed subgroup of $\GL(n;\ol{k})$. Associated to $\B{G}$ is its \emph{ideal of vanishing} $\Fr{a}$, and we say $\B{G}$ is \emph{$L$--defined} or \emph{$L$--algebraic} if $\Fr{a}$ is generated over a subfield $L$ of $\ol{k}$. For any subring $R$ of $\ol{k}$ whose field of fractions is $L$, the subgroup $\B{G}(R) = \B{G} \cap \GL(n;R)$ is well-defined up to commensurability and called the \emph{$R$--points} of $\B{G}$. We say $\B{G}$ is a \emph{real algebraic group} if $\B{G}$ is the $\R$--points of an algebraic group. 

\begin{lemma}[\cite{Raghunathan72}]\label{L:323}
Let $\Ga<\B{G}(k)$ and $\B{G}(\Cal{O}_k)$ be commensurable and 
\[ f\co \B{G} \lra \B{G}^\prime \] 
be a $k$--homomorphism of $k$--algebraic groups. Then there exists $\Ga^\prime<\B{G}^\prime(k)$ commensurable with $\B{G}^\prime(\Cal{O}_k)$ such that $f(\Ga)$ is contained in $\Ga^\prime$.
\end{lemma}

For a number field $k$, a linear $k$--algebraic group $\B{G}$, and an ideal $\Fr{b}$ of $\Cal{O}_k$, the kernel of the \emph{reduction homomorphism} 
\[ r_{\Fr{b}}\co \B{G}(\Cal{O}_k) \lra \B{G}(\Cal{O}_k/\Fr{b}) \] 
is a finite index subgroup called a \emph{principal congruence subgroup}. Any subgroup of $\B{G}(\Cal{O}_k)$ containing a principal congruence subgroup will be called a \emph{congruence subgroup}.

\subsection{Corestriction}

Let $k$ be a number field and $k_{gal}$ its Galois closure over $\Q$. For a $k$--algebraic group $\B{G}$ and finite $k$--generating set $P_1(T),\dots,P_r(T)$ of $\Fr{a}$, each element $\si$ in $\Gal(k_{gal}/\Q)$ yields a new set of $k$--polynomials $^\si P_1(T),\dots,~^\si P_r(T)$. The ideal $^\si\Fr{a}$ generated by these polynomials is the ideal of vanish of a $k$--algebraic group $^\si \B{G}$, and for any two automorphisms $\si_1,\si_2 \in \Gal(k_{gal}/\Q)$ equivalent modulo $\Gal(k_{gal}/k)$, the groups $^{\si_1}\B{G},~^{\si_2}\B{G}$ are $k$--isomorphic. We define the group
\[ \res_{k/\Q}(\B{G}) = \prod_{\si \in \Gal(k_{gal}/\Q)/\Gal(k_{gal}/k)} ~^{\si}\B{G}. \]
By construction $\res_{k/\Q}(\B{G})$ is invariant under the action of $\Gal(k_{gal}/\Q)$, and so $\res_{k/\Q}(\B{G})$ is $\Q$--algebraic and the groups $\res_{k/\Q}(\B{G})(\Z)$ and $\res_{k/\Q}(\B{G}(\Cal{O}_k))$ are commensurable. We call this process \emph{restriction of scalars} or \emph{corestriction}.

\subsection{Arithmetic lattices}

For a $H$ locally compact group with right Haar measure $\mu$, a \emph{lattice} $\La$ of $H$ is a discrete subgroup such that $\mu$ descends to a finite measure on $H/\La$. When $H/\La$ is compact, $\La$ is called \emph{cocompact} and otherwise \emph{noncocompact}. It is fundamental result of Borel and Harish-Chandra \cite{BorelChandra62} that $\B{G}(\Z)$ is a lattice in $\B{G}(\R)$ for any $\Q$--algebraic semisimple group $\B{G}$. More generally, for a real Lie group $G$, we say $\La<G$ is an \emph{arithmetic lattice} if there exists a semisimple $\Q$--algebraic group $\B{G}$, compact Lie groups $K_1$ and $K_2$, and an exact sequence
\[ \xymatrix{ 1 \ar[r] & K_1 \ar[r] & \B{G}(\R) \ar[r]^\psi & G \ar[r] & K_2 \ar[r] & 1 } \]
such that $\psi(\B{G}(\Z))$ is commensurable with $\La$.

\subsection{$k$--forms of a Lie group}

A \emph{$k$--form} of  a Lie group $G$ is a real $k$--algebraic group $\B{G}$ and a Lie epimorphism $\rho\co \B{G} \lra G$ with compact kernel. Applying restriction of scalars to $\B{G}$, we obtain a real $\Q$--algebraic group $\res_{k/\Q}(\B{G})$ and the diagram:
\[ \xymatrix{ & & & 1 \ar[d] & \\ & & & \ker \rho \ar[d] & \\ 1 \ar[r] & \prod_{\si \ne \id} ~^{\si}\B{G} \ar[r] & \res_{k/\Q}(\B{G}) \ar[r]^\pi \ar@{-->}[rd]_{\rho \circ \pi} & \B{G} \ar[r] \ar[d]^\rho & 1 \\ & & & G \ar[d] & \\ & & & 1 & } \]
When $\ker \pi$ is compact and $\B{G}$ is semisimple, $\rho(\B{G}(\Cal{O}_k))$ is an arithmetic lattice in $G$ and we call such a $k$--form of $G$ \emph{admissible}.

\section{Arithmetic lattices in the classical rank one groups}

In our investigation of cusp cross-sections of arithmetic lattices, we require a classification theorem for noncocompact arithmetic lattices in $\PO(n,1)$, $\PU(n,1)$, and $\PSp(n,1)$. This is summarized in the following theorem often attributed to Weil.

\begin{thm}\label{T:ClassicationNCL}
Let $\La$ be a noncocompact arithmetic lattice in $\Isom(\Hy_X^n)$.
\begin{description}
\item[(a)] If $X=\R$, then $\La$ is commensurable with $\PO(B;\Z)$, where $B$ is a signature $(n,1)$ bilinear form defined over $\Q$.
\item[(b)] If $X=\C$, then $\La$ is commensurable with $\PU(H;\Cal{O}_k)$, where $H$ is a hermitian form of signature $(n,1)$ defined over an imaginary quadratic number field $k$.
\item[(c)] If $X=\BB{H}$, then $\La$ is commensurable with $\PU(H;\Cal{O})$, where $H$ is a hermitian form of signature $(n,1)$ defined over a ramified quaternion $\Q$--algebra $A$ and $\Cal{O}$ is a $\Z$--order in $A$. 
\end{description}
\end{thm}

\begin{rem} The converse to \refT{ClassicationNCL} holds modulo a few exceptional cases. Namely, the groups $\PO(B;\Z)$, $\PU(H;\Cal{O}_k)$, and $\PU(H;\Cal{O})$ are noncompact lattices except possibly for $n=2,3$ and $X=\R$ when some of these lattices are cocompact.
\end{rem}

We prove \refT{ClassicationNCL} (b) as we could not locate a proof in the literature. For (a) and (c), the reader is directed to \cite{WitteA}---for (a), see also \cite{LiMillson93}. 

\subsection{Involutions and matrices over algebras}

For a central simple $E$--algebra $B$, an \emph{involution} on $B$ is an order two map $*\co B \lra B$ such that $(x+y)^* = x^* + y^*$ and $(xy)^* = y^*x^*$. Given an involution $*$ on an algebra $B$, we obtain an involution $\star$ on $\Mat(n;B)$ given by $(a_{ij})^\star = (a_{ji}^*)$, and call this \emph{$*$--conjugate transpose}. An element $b \in \GL(n;B)$ is \emph{$\star$--symmetric} if $b^\star=b$, and associated to $b$ is the group
\[ \B{G}_b(B) = \set{x \in \SL(n;B)~:~b^{-1}x^\star b x = 1}. \] 

\subsection{Arithmetic lattices in $\SU(n,1)$}

Let $E/F$ be a CM field, $(r,d) \in \N^2$ with $rd=n+1$, and $B$ a cyclic $E$--division algebra $B$ of degree $d$ that admits an involution $*$ whose restriction to $E$ is the nontrivial Galois involution $\te \in \Gal(E/F)$. For a $\star$--symmetric element $b$ in $\Mat(r;B)$ with associated group $\B{G}_b(B)$ and each embedding $\tau$ of $E$, we obtain a new group $^\tau\B{G}_b(~^\tau B)$. By Weil's theorem on involutions \cite{Weil60}, the group $^\tau\B{G}_b(~^\tau B \otimes \C)$, which we denote by $^\tau\B{G}_b(\R)$, is isomorphic as a real algebraic group to $\SU(r_\tau,s_\tau)$ for some $r_\tau,s_\tau$ such that $r_\tau+s_\tau=n+1$. We say $(B,b)$ is \emph{admissible} if
\[ ^{\tau_j}\B{G}_b(\R) = \begin{cases} \SU(n,1), & \tau_j=\tau_1 \\ \SU(n+1), & \text{otherwise}. \end{cases} \]
For any admissible $(B,b)$ and $\Cal{O}_E$--order $\Cal{O}$ of $B$, the group $\B{G}_b(\Cal{O})$ is a lattice in $\SU(n,1)$ via the inclusion provided by any isomorphism between $\B{G}_b(\R)$ and $\SU(n,1)$. These lattices are arithmetic and according to Tits \cite[p. 33--66]{BorelMostow66}, up to wide commensurability, comprise the full class of arithmetic lattices.

\subsection{The proof of \refT{ClassicationNCL}}

The proof of \refT{ClassicationNCL} requires the following theorem (see \cite{BorelChandra62}, \cite{MostowTamagawa62}).  

\begin{thm}[Godement's compactness criterion]\label{T:Godement}
\begin{description}
\item[(a)]
Let $\B{G}$ be a semisimple $\Q$--algebraic group and $\La$ an arithmetic lattice in $\B{G}$. Then $\La$ is cocompact if and only if $\La$ contains no non-trivial unipotent elements.
\item[(b)]
Let $G$ be a noncompact semisimple Lie group with an admissible $k$--form $\B{G}$. If $\B{G}(\Cal{O}_k)$ is noncocompact, then $k=\Q$.
\end{description}
\end{thm}

\begin{pf}[Proof of \refT{ClassicationNCL}] It suffices to prove this for lattices in $\SU(n,1)$. In addition, since noncocompactness is an invariant of the wide commensurability class, it suffices to determine when the groups $\B{G}_b(\Cal{O})$ are noncocompact. To this end, let $(B,b)$ be admissible for the pair $(r,d)$ over the CM field $E/F$ such that $\B{G}_b(\Cal{O})$ is noncocompact, and note our goal is to show $F=\Q$ and $d=1$. The latter is achieved by part (b) of \refT{Godement}, and the admissibility assumption in this case is simply that $\B{G}_b(\R)$ be isomorphic to $\SU(n,1)$.
 
By definition, there exists a cyclic extension $L/E$ of degree $d$ such that $B=(L/E,\te,\al)$. The field $L$ is totally imaginary, $L/K$ is a CM field for a unique subfield $K$ of $L$, and we have the field diagram
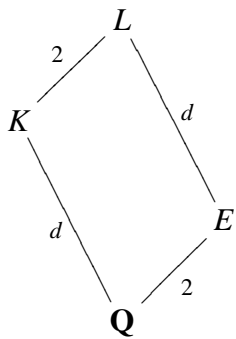
\[ \xymatrix{ & L \ar@{-}[ld]_2 \ar@{-}[rdd]^d & \\ K \ar@{-}[rdd]_d & & \\ & & E \ar@{-}[ld]^2 \\ & \Q & } \] 
Up to wide commensurability, we may assume (\cite{WitteA})
\[ b = \begin{pmatrix} \al_1 & 0 & \dots & 0 \\ 0 & \al_2 & \dots & 0 \\ \vdots & \vdots & \ddots & \vdots \\ 0 & 0 & \dots & \al_r \end{pmatrix} \]
with $\al_1,\al_2,\dots,\al_r \in K$. By \refT{Godement} (a), $\B{G}_b(\Cal{O})$ contains a nontrivial unipotent element, and so by a change of $B$--basis, we may further assume $b$ has the form (\cite{WitteA})
\[ b = \begin{pmatrix} \be_1 & 0 & \dots & 0 & 0 & 0 \\ 0 & \be_2 & \dots & 0 & 0 & 0 \\ \vdots & \vdots & \ddots & \vdots & \vdots & \vdots \\ 0 & 0 & \dots & \be_{r-2} & 0 & 0 \\ 0 & 0 & \dots & 0 & 1 & 0 \\ 0 & 0 & \dots & 0 & 0 & -1 \end{pmatrix} \]  
with $\be_1,\be_2,\dots,\be_{r-2} \in K$. The group $\B{G}_b(\R)$ is given by extending the coefficients of $B$ from $E$ to $\C$, and $\B{G}_b(\R)$ is $\SU(H)$, where $H$ is the image of $b$ under the embedding $\Mat(r;B) \lra \Mat(r;B\otimes_E \C)$. The image of $b$ is given explicitly by
\[ \bop _{j=1}^{r-2}\begin{pmatrix} \be_1 & 0 & \dots & 0 \\ 0 & \te(\be_1) & \dots & 0 \\ \vdots & \vdots & \ddots & \vdots \\ 0 & 0 & \dots & \te^{d-1}(\be_1) \end{pmatrix} \op \bop_{j=1}^d \begin{pmatrix} 1 & 0 \\ 0 & -1 \end{pmatrix}, \]
where by this we mean the block diagonal matrix with these blocks. According to the admissibility assumption, $H$ must have signature $(n,1)$ while the signature of the image of $b$ is $(p+d,q+d)$, for some nonnegative integers $p,q$. In particular, this can happen if and only if $d=1$ and $q=0$, and the resulting groups are precisely those in the statement of (b). 
\end{pf} 

\section{Solvable groups and their geometries}

\subsection{Bieberbach groups}

We denote the affine, Euclidean, and similarity groups of $\R^n$ by $\Aff(n)$, $\Euc(n)$, and $\Sim(n)$. Given a discrete, torsion free subgroup $\Ga$ of $\Euc(n)$, the quotient $\R^n/\Ga$ is a smooth manifold, and those $\Ga$ whose quotient manifold is compact will be called \emph{Bieberbach groups}. We denote the space of faithful representations of $\Ga$ in $\Aff(n)$ with Bieberbach images by $\Cal{R}_f(\Ga)$ (i.e., an $\Aff(n)$--conjugate of a Bieberbach group). Finally, $\Cal{F}_f(\Ga), \Cal{S}_f(\Ga)$ will denote the subspaces of the $\Euc(n)$ and $\Sim(n)$--character spaces consisting of the faithful characters whose image is Bieberbach.

Associated to each maximal compact subgroup $K$ of $\GL(n;\R)$ is the \emph{orthogonal affine group} $\Ort_K(n)= \R^n \rtimes K$. As each $K$ is conjugate in $\GL(n;\R)$ to $\Ort(n)$, $K$ is $\Ort(B_K)$ for some symmetric, positive definite, bilinear form $B_K$ on $\R^n$. When $B_K$ is $\Q$--defined (i.e. $\Q$--valued on some $\R$--basis), $\Ort_K(n)$ is $\Q$--defined and we call subgroups commensurable with $\Ort_K(n;\Z)$ \emph{$\Q$--arithmetic}. We say $\rho$ in $\Cal{R}_f(\Ga)$ is \emph{$\Q$--arithmetic} if there exists a $\Q$--defined orthogonal affine group $\Ort_K(n)$ such that $\rho(\Ga)$ is a $\Q$--arithmetic subgroup of $\Ort_K(n)$, and denote the subspace of $\Cal{R}_f(\Ga)$ of $\Q$--arithmetic representations by $\Cal{R}(\Ga;\Q)$. 
 
\subsection{Almost Bieberbach groups}

Let $X$ be $\R,\C$, or $\BB{H}$, $\ell = \dim_\R X$, $\innp{z,w}$ the standard inner product on $X^n$, and $\om(z,w) = 2\Ima\innp{z,w}$. The \emph{$X$--Heisenberg group} $\Fr{N}_{\ell n-1}(X)$ is the extension
\[ 1 \lra \Ima X \lra \Fr{N}_{\ell n-1}(X) \lra X^{n-1} \lra 1 \]
with trivial holonomy and associated 2--cocycle $\om$. When $X=\C$, we call this group the \emph{$(2n-1)$--dimensional Heisenberg group} and for $X=\BB{H}$, we call this group the \emph{$(4n-1)$--dimensional quaternionic Heisenberg group}.

The automorphism group of the $X$--Heisenberg group $\Aut(\Fr{N}_{\ell n-1}(X))$ splits as $\Inn(\Fr{N}_{\ell n-1}(X)) \rtimes
\Out(\Fr{N}_{\ell n-1}(X))$, and the $\om$--nondegenerate vectors in $X^{n-1}$ are in bijection with the nontrivial inner automorphisms---this is the whole of $X^{n-1}\smin \set{0}$ except when $X=\R$. The outer automorphism group is comprised of three types of automorphisms. The first type of automorphism is a \emph{symplectic rotation} given by 
\[ S(\xi,t) = (S\xi,t) \] 
for $S \in \Sp(\om)$.  The second type of automorphism is a \emph{Heisenberg dilation} given by
\[ d(\xi,t) = (d\xi,d^2t) \] 
for $d \in \R^\times$. Finally, we have \emph{$X$--scalar conjugation} given by 
\[ \zeta(\xi,t) = (\zeta^{-1}\xi\zeta, \zeta^{-1}t\zeta) \] 
for $\zeta \in X^\times$. The outer automorphism group is generated by these three automorphisms, and in total
\[ \Out(\Fr{N}_{\ell n-1}(X)) = \begin{cases} \GL(n-1;\R), & X=\R \\ \Sp(2n-2) \times \R^\times, & X=\C \\ \Sp(\om) \times \R^\times \times \BB{H}^\times, & X=\BB{H}. \end{cases} \]  

The maximal compact subgroups of $\Aut(\Fr{N}_{\ell n -1}(X))$ are of the form
\[ M(X) = \begin{cases} \Ort(B_{M(X)}), & X=\R \\ \innp{\Uni(H_{M(X)}),\iota}, & X=\C \\ \Uni(H_{M(X)}) \times S, & X=\BB{H}, \end{cases} \]
where $B_{M(X)}$ is a symmetric, positive definite bilinear form, $H_{M(X)}$ is a positive definite hermitian form with $\Ima H_{M(X)}
= \om$, and $S$ is the unit sphere in $\BB{H}$ (equipped possibly with a nonstandard quaternionic structure). As all maximal compact subgroups are conjugate in $\Aut(\Fr{N}_{\ell n -1}(X)$, each $M(X)$ is Lie isomorphic to
\[ M_s(X) = \begin{cases} \Ort(n-1), & X=\R \\ \innp{\Uni(n-1),\iota}, & X=\C \\ \Sp(n-1) \times \SO(3), & X=\BB{H}. \end{cases} \]

For a given maximal compact subgroup $M$, we call the group 
\[ \Fr{N}_{\ell n-1}(X) \rtimes M \] 
a \emph{unitary affine group} and denote this group by $U_M(n-1;X)$. We call the group 
\[ \Fr{N}_{\ell n-1}(X) \rtimes (M(X) \times \R^+) \]
an \emph{$X$--Heisenberg similarity group} and denote this group by $S_M(n-1;X)$. Finally, we call the group 
\[ \Fr{N}_{\ell n-1}(X) \rtimes \Aut(\Fr{N}_{\ell n-1}(X)) \]
the \emph{$X$--Heisenberg affine group} and denote it by $\Aff(\Fr{N}_{\ell n-1}(X))$.

By an \emph{almost Bieberbach group} (or \emph{AB-group} for short) modelled on $\Fr{N}_{\ell n-1}(X)$, we mean a discrete, torsion free subgroup $\Ga$ of $\Aff(\Fr{N}_{\ell n-1}(X))$ such that the quotient $\Fr{N}_{\ell n-1}(X)/\Ga$ is compact and $\Ga \cap \Fr{N}_{\ell n-1}(X)$ is a finite index subgroup of $\Ga$. Every AB-group modelled on $\Fr{N}_{\ell n-1}(X)$ is determined by the short exact sequence 
\[ 1 \lra L \lra \Ga \lra \te \lra 1, \]
where $L = \Ga \cap \Fr{N}_{\ell n-1}(X)$ and $\abs{\te}<\iny$. We call $L$ the \emph{Fitting subgroup of $\Ga$} and $\te$ the \emph{holonomy group of $\Ga$}. The above exact sequence induces an injective homomorphism $\vp\co \te \lra \Out(\Fr{N}_{\ell n-1}(X))$ called the \emph{holonomy representation of $\te$}. Since $\te$ is finite, this is conjugate into a representation $\vp\co \te \lra M(X)$ for any $M(X)$ and yields a faithful representation $\rho\co \Ga \lra U_M(n-1;X)$ for any $M(X)$.

\subsection{Flat, almost flat, and infrasolv manifolds}

We say a connected, closed manifold $M$ is \emph{flat} if $M$ is diffeomorphic to $\R^n/\Ga$ for some Bieberbach group. The flat metric $g$ induced by the standard inner product $\innp{\cdot,\cdot}$ on $\R$ supplies $M$ with a flat metric called the \emph{associated flat structure}, and $\Cal{F}(M)$ will denote the space of all isometry classes of such metrics. We say two flat structures $g_1$ and $g_2$ are \emph{similar} if there exists an isometry $f\co (M,g_1) \lra (M,\al g_2)$, for some $\al \in \R^+$. We denote the similarity class of a flat metric $g$ by $[g]$ and the space of all similarity classes by $\Cal{S}(M)$.

\begin{thm}[\cite{Thurston97}] The space of flat isometry classes on $M$ is $\Cal{F}_f(\pi_1(M))$. The space of flat similarity classes on $M$ is $\Cal{S}_f(\pi_1(M))$.
\end{thm}

Any faithful representation $\rho$ of $\pi_1(M)$ into $\Ort_K(n)$ whose image is Bieberbach endows $M$ with a flat metric induced from the form $B_K$.

For a complete Riemannian manifold $(M^n,g)$, let $d(g)$, $c^-(g)$ and $c^+(g)$ denote the diameter of $M$ and the lower and
upper bounds of the sectional curvature of $M$, respectively, and set $c(g)$ to be the maximum of $\abs{c^+}$ and $\abs{c^-}$. We say $M$ is \emph{almost flat} if there exists a family of complete Riemannian metrics $g_j$ on $M$ such that 
\[ \lim_{j \lra \iny} d(g_j)^2c(g_j) = 0. \]
Gromov \cite{Gromov78} proved every compact almost flat manifold is of the form $N/\Ga$, where $N$ is a connected, simply connected nilpotent Lie group and $\Ga$ is an AB-group modelled on $N$. In the sequel we refer to compact almost flat manifolds as \emph{infranil manifolds} modelled on $N$, where $N$ is the connected, simply connected nilpotent cover. In the event the fundamental group is a lattice in $N$, we call such manifolds \emph{nil manifolds} modelled on $N$.

For a simply connected, connected, solvable, Lie group $S$, a discrete torsion free subgroup $\Ga$ of $\Aff(S)$ is an \emph{infrasolv group modelled on $S$} if $\Ga \cap S$ is finite index in $\Ga$ and $S/\Ga$ is a compact manifold. A smooth manifold diffeomorphic with $S/\Ga$ for some infrasolv group will be called an \emph{infrasolv manifold modelled on $S$}, and we require the following rigidity result of Mostow \cite{Mostow54} for these manifolds.

\begin{thm}\label{T:Rig} If $M_1$, $M_2$ are infrasolv manifolds with $\pi_1(M_1) \cong  \pi_1(M_2)$, then $M_1$ is diffeomorphic to $M_2$.
\end{thm}

\subsection{$\Nil$ and $\Sol$ geometry}

The 3--dimensional geometries $\Nil$ and $\Sol$ play a prominent role in this thesis, and for completeness, a short introduction is provided here.

The Heisenberg group $\Fr{N}_3$ can also be viewed as the subgroup of $\SL(3;\R)$ of matrices of the form
\[ \begin{pmatrix} 1 & x & t \\ 0 & 1 & y \\ 0 & 0 & 1 \end{pmatrix}. \]
Identifying $S^1$ with the rotations in the $xy$--plane, an \emph{orientable $\Nil$ 3--manifold} is a manifold of the form
$\Fr{N}_3/\Ga$, where $\Ga$ is a discrete subgroup of $\Fr{N}_3 \rtimes S^1$ which acts freely. As we will have need for this in the sequel, we must also consider an orientation reversing involution given by
\[ \wt{\iota}\begin{pmatrix} 1 & 0 & 0 \\ 0 & 1 & y \\ 0 & 0 & 1 \end{pmatrix} = \begin{pmatrix} 1 & 0 & 0 \\ 0 & 1 & -y \\ 0 & 0 & 1
\end{pmatrix}, \quad \wt{\iota} \begin{pmatrix} 1 & x & 0 \\ 0 & 1 & 0 \\ 0 & 0 & 1 \end{pmatrix} = \begin{pmatrix} 1 & x & 0 \\ 0 & 1 & 0 \\
0 & 0 & 1 \end{pmatrix}. \]
As automorphisms of lattices of $\Fr{N}_3$ uniquely determine automorphisms of $\Fr{N}_3$ by Mal'cev rigidity, this determines a continuous isomorphism of $\Fr{N}_3$.  A \emph{$\Nil$ 3--manifold} is a manifold of the form $\Fr{N}_3/\Ga$, where $\Ga$ is a discrete subgroup of $\innp{\Fr{N}_3\rtimes S^1,\wt{\iota}}$ which acts freely.

The solvable Lie group $\Sol$ is defined to be topologically $\R^2 \times \R^+$ with group operation defined by
\[ (x_1,y_1,t_1)\cdot(x_2,y_2,t_2) \bdef (x_1 + e^{t_1}x_2,y_1+  e^{-t_1}y_2, t_1+t_2). \]
By a \emph{$\Sol$ 3--manifold}, we mean a manifold $M$ diffeomorphic with $\Sol/\Ga$, where $\Ga$ is a discrete, torsion free subgroup of $\Aff(\Sol)$ such that $\Sol/\Ga$ is compact and $[\Ga:\Ga \cap \Sol]<\iny$. 

\section{$X$--Hyperbolic geometry}

For $X=\R,\C$, or $\BB{H}$, the classical groups $\SO(n,1)$, $\SU(n,1)$, and $\Sp(n,1)$ produce the symmetric spaces $\Hy_\R^n$, $\Hy_\C^n$, and $\Hy_\BB{H}^n$, known collectively as $X$--hyperbolic $n$--space. For an explicit description, we equip $X^{n+1}$ with a hermitian form $H$ of signature $(n,1)$. \emph{$X$--hyperbolic $n$--space} is the (left) $X$--projectivization of the $H$--negative vectors endowed with the Bergman metric associated to $H$. We denote $X$--hyperbolic $n$--space together with this metric by $\Hy_X^n$ and say $\Hy_X^n$ is \emph{modelled on $H$} and call $H$ a \emph{model form}. The \emph{boundary} of $\Hy_X^n$ in $PX^{n+1}$ is the $X$--projectivization of the $H$--null vectors. We denote this set by $\prt \Hy_X^n$, which is topologically just $S^{\ell n}$---see \cite[p. 265]{BridsonHaefliger99}.

The spaces constructed in this way yield every symmetric space of real rank--1 except for the exceptional \emph{Cayley hyperbolic plane} $\Hy_\B{O}^2$. We shall only make use of the fact that $\Isom(\Hy_\B{O}^2)$ has a faithful linear representation and refer the reader to \cite{Allcock99} for more on the Cayley hyperbolic plane.

The isometry group of $\Hy_X^n$ is denoted by $\Isom(\Hy_X^n)$, and in each setting is locally isomorphic to $\SU(H)$. More precisely, 
\[ \Isom(\Hy_X^n) = \begin{cases} \innp{\PU(H)_0,\iota}, & X=\R,\C \\ \PU(H), & X=\BB{H}, \end{cases} \]
where $\iota$ is an involution induced by an inversion in the real case and complex conjugation in the complex case.

\subsection{The Iwasawa decomposition of $\Isom(\Hy_X^n)$}

The isometry group of $X$--hyperbolic $n$--space decomposes as $KAN$ via the \emph{Iwasawa decomposition} (see \cite[p. 311--313]{BridsonHaefliger99}). The factor $N$ is isomorphic to the $X$--Heisenberg group $\Fr{N}_{\ell n-1}(X)$ and all isomorphisms arise in the following fashion. Let $H$ be a model hermitian form for $X$--hyperbolic $n$--space and $V_\iny$ be the $H$--orthogonal complement of a pair of $X$--linearly independent $H$--null vectors $v_0$ and $v_\iny$ in $X^{n+1}$. For any maximal compact $M(X)$ of $\Aut(\Fr{N}_{\ell n - 1})$ and associated positive definite hermitian form $H_{M(X)}$, let 
\[ \psi\co (X^{n-1},H_{M(X)}) \lra (V_\iny,H_{|V_\iny}) \] 
be any isometric $X$--isomorphism. This induces a map  $\eta\co X^{n-1} \lra N$ defined by 
\[ \eta(\xi) = \exp(\psi(\xi) v_\iny^* - v_\iny\psi(\xi)^*), \] 
where $xy^*(\cdot) = H(\cdot,y)x$ is the hermitian outer pairing of $x$ and $y$ with respect to the hermitian form $H$. This extends to all of $\Fr{N}_{\ell n-1}(X)$ as these elements generate $\Fr{N}_{\ell n-1}(X)$. In fact, this extends to $\eta\co S_M(n-1;X) \lra \Isom(\Hy_X^n)$, and produces the equality $\eta(S_M(n-1;X)) = \Stab(v_\iny)$. 

\subsection{The upper half plane model}

Viewing $\Hy_X^n$ as the coset space $\Isom(\Hy_X^n)/K$, where $K$ is a maximal compact subgroup, topologically $\Hy_X^n$ is $A \times N$, and we call this the \emph{upper half plane} model for $\Hy_X^n$. As $A=\R^+$, $\Hy_X^n$ has a foliation 
\[ \bu_{t \in \R^+} \set{t} \times N \] 
whose leaves are called \emph{horospheres} and are said to be \emph{centered} at $v_\iny$ when $N$ arises from $v_\iny$.

\begin{figure}[h!]
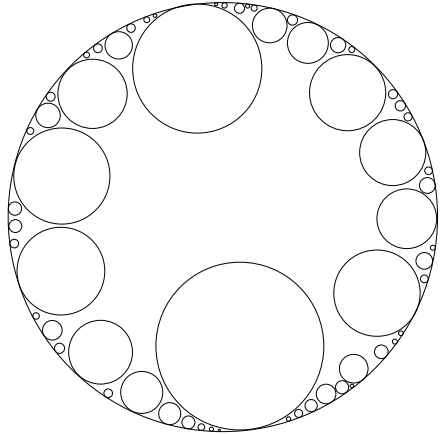

\begin{center}
\ig{NeuteredSpace}
\caption{Some horospheres in the disk model for the hyperbolic plane.}
\end{center}
\end{figure}  

\subsection{The Kazhdan--Margulis theorem and cusp cross-sections}

For a lattice $\La$ in $\Isom(\Hy_X^n)$ with associated orbifold $M=\Hy_X^n/\La$, we say $M$ (or $\La$) has a \emph{cusp} at $v \in \prt \Hy_X^n$ if $\La \cap N \ne \set{1}$ and $\set{t} \times N$ is centered at $v$.

For a lattice $\La$ in $\Isom(\Hy_X^n)$ with cusp at $v$, we define the \emph{maximal peripheral subgroup of $\La$ at $v$} to be the subgroup 
\[ \tri_v(\La) = \Stab(v) \cap \La. \] 
By the Kazhdan--Margulis theorem \cite{KazhdanMargulis68}, $\tri_v(\La)$ is virtually nilpotent with the maximal, torsion free, nilpotent subgroup $\tri_v(\La) \cap N$, and there exists a horosphere $\set{t} \times N$ such that $(\set{t}\times N)/\tri_v(\La)$ is embedded in $\Hy_X^n/\La$. We call $(\set{t} \times N)/\tri_v(\La)$ a \emph{cusp cross-section} of the cusp at $v$. 

More generally, for any $v \in \prt \Hy_X^n$, there are three possibilities for $\tri_v(\La)$: 
\begin{description}
\item[(1)] $\tri_v(\La)$ is finite.
\item[(2)] $\tri_v(\La)$ is virtually cyclic with cyclic subgroup generated by a loxodromic isometry. 
\item[(3)] $\tri_v(\La)$ is an AB-group modelled on $\Fr{N}_{\ell n-1}(X)$.
\end{description}
In the sequel, we refer to this trichotomy as the \emph{stabilizer trichotomy}. 

\chapter{Separable subgroups and geometric applications}

In this chapter, we prove the subgroup separability results from the introduction. This task begins the following basic lemma---for a general reference, we refer the reader to \cite{LongReid05}.

\begin{lemma}\label{L:Basic}
\begin{description}
\item[(a)]
Let $H<K<G$ and $H$ be separable in $G$ with $[K:H]<\iny$. Then $K$ is separable in $G$.
\item[(b)]
Let $H,L$ be separable subgroups of $G$. Then $H \cap L$ is separable in $G$.
\item[(c)]
Let $H,G_0$ be subgroups of $G$ with $[G:G_0]<\iny$. Then $H$ is separable in $G$ if and only if $(G_0 \cap H)$ is separable in $G_0$.
\item[(d)]
Let $H<G_0<G$ and $H$ be separable in $G$. Then $H$ is separable in $G_0$. 
\end{description}
\end{lemma}

\begin{pf} For brevity, we only prove (b) and (c), as the proofs of (a) and (d) are of a similar flavor. 

Part (b). Let $\ga \in G \smin (H \cap L)$ and assume $\ga \notin H$. Since $H$ is separable in $G$, there exists a finite index subgroup $K$ of $G$ with $H$ contained in $K$ and $\ga \notin K$. Visibly, $K$ separates $\ga$ and $H \cap L$. For the alternative $\ga \notin L$, an identical argument is made using $L$.

Part (c). The direct implication is immediate by (b), since $H \cap G_0$ is separable in the larger group $G$. For the reverse implication, to show $H$ is separable in $G$, by (a) it suffices to show $G_0 \cap H$ is separable in $G$. For $g \in G \smin (G_0 \cap H)$, there are two cases to consider. If $g \notin G_0$, then $G_0$ separates $G_0 \cap H$ and $g$. Otherwise, if $g \in G_0$, since $G_0 \cap H$ is separable in $G_0$, there exists a finite index subgroup $K$ of $G_0$ separating $G_0 \cap H$ and $g$ in $G_0$. As $[G:G_0]<\iny$, $K$ is also a finite index subgroup of $G$ and thus separates $G_0 \cap H$ and $g$ in $G$.
\end{pf} 

\section{Line stabilizers and Chevalley's theorem}

In this section, we prove the following result which bears \refT{Borel} as a corollary.

\begin{thm}[Closed stabilizer separability]\label{T:CSS}
Let $k$ be a number field, $\B{G}$ a linear $k$--algebraic group, and $\B{H}$ a closed $k$--algebraic subgroup such that every finite index subgroup of $\B{H}(\Cal{O}_k)$ is congruence. Then a subgroup of $\B{H}(\Cal{O}_k)$ is separable in $\B{H}(\Cal{O}_k)$ if and only if it is separable in $\B{G}(\Cal{O}_k)$.
\end{thm}

To prove this we employ the following theorem of Chevalley---see \cite{PlatonovRapinchuk94}.

\begin{thmnl}[Chevalley's theorem]
Let $\B{G}$ be a linear $k$--algebraic group and $\B{H}$ a closed $k$--algebraic subgroup. Then there exists a faithful $k$--homomorphism
$\vp\co \B{G} \lra \GL(m;\ol{k})$ and a $k$--defined line $\ell$ in $\ol{k}^m$ such that $\Stab_{\B{G}}(\ell) = \B{H}$.
\end{thmnl}

\begin{pf}[Proof of \refT{CSS}] Let $H$ be a subgroup of $\B{H}(\Cal{O}_k)$. By \refL{Basic} (d), if $H$ is separable in $\B{G}(\Cal{O}_k)$, then $H$ is separable in $\B{H}(\Cal{O}_k)$. For the converse, assume $H$ is separable in $\B{H}(\Cal{O}_k)$. According to Chevalley's theorem, there exists a faithful $k$--homomorphism $\vp\co \B{G} \lra \GL(m;\ol{k})$ and a $k$--defined line $\ell$ in $\ol{k}^m$ such that $\Stab_{\B{G}}(\ell) = \B{H}$. We can select a basis $v_1,\dots,v_m$ of $\ol{k}^m$ such that both $\ell$ is spanned by $v_m$ and the inner automorphism $\tau$ of $\GL(m;\ol{k})$ induced by the change of basis from the standard basis $\set{e_1,\dots,e_m}$ to $\set{v_1,\dots,v_m}$ is $k$--defined. Let $\vp_\tau$ denote the $k$--defined homomorphism of $\B{G}$ into $\GL(m;\ol{k})$ given by $\tau \circ \vp$. By \refL{323}, $\vp_\tau^{-1}(\vp_\tau(\B{G}(\Cal{O}_k)) \cap \GL(m;\Cal{O}_k))$ is a finite index subgroup of $\B{G}(\Cal{O}_k)$ and we denote this subgroup by $G_0$. According to \refL{Basic} (c), it suffices to separate $H_0=G_0 \cap H$ in $G_0$, and this is achieved as follows. For $\ga \in G_0\smin H_0$, we have two cases to consider. First, if $\vp(\ga)$ is not in $\Stab(\ell)$, then $\vp_\tau(\ga)_{i,m} \ne 0$ for some $1 \leq i < m$, where $\vp_\tau(\ga)_{i,m}$ denotes the $(i,m)$--coefficient of the matrix $\vp_\tau(\ga)$. As $\Cal{O}_k$ is Dedekind, there are only finitely many prime ideals $\Fr{p}_1,\dots,\Fr{p}_{j_\ga}$ of $\Cal{O}_k$ such that $\la_{i,m} = 0 \mod \Fr{p}_j$. In particular, so long as $\Fr{p}$ is an ideal different from $\Fr{p}_1,\dots,\Fr{p}_{j_\ga}$, the number $r_\Fr{p}(\vp_\tau(\ga))_{i,m}$ is nonzero. In contrast, since $\vp_\tau(H_0)$ stabilizes the line $\ell$, every element $h$ in $H_0$ has the property that $r_\Fr{p}(\vp_\tau(h))_{i,m}=0$. Therefore, $r_\Fr{p}(\ga) \notin r_\Fr{p}(H_0)$ and thus $r_\Fr{p}^{-1}(r_\Fr{p}(H_0)) \cap G_0$ is a finite index subgroup of $G_0$ separating $H_0$ and $\ga$. This leaves the alternative $\vp(\ga) \in \Stab(\ell)$ and the remaining task of separating $H_0$ from $\ga$ in $G_0$. By assumption $H$ is separable in $\B{H}(\Cal{O}_k)$, and according to \refL{Basic} (c), $H \cap G_0$ is separable in $\B{H}(\Cal{O}_k) \cap G_0$. As $H_0=H\cap G_0$, $H_0$ is separable in $\B{H}(\Cal{O}_k) \cap G_0$ and hence there must exist a finite index subgroup $K_0$ of $\B{H}(\Cal{O}_k) \cap G_0$ separating $H_0$ and $\ga$. Since $K_0$ is a finite index subgroup of $\B{H}(\Cal{O}_k)$, by \refT{Ch}, there exists an ideal $\Fr{a}$ in $\Cal{O}_k$ such that $\ker r_{\B{H},\Fr{a}} < K_0$, where $r_{\B{H},\Fr{a}}$ is the modulo $\Fr{a}$ reduction for $\B{H}(\Cal{O}_k)$. It is a simple matter to see that $r_{\B{H},\Fr{a}}(\ga) \notin r_{\B{H},\Fr{a}}(K_0)$, and hence the desired finite index subgroup of $G_0$ separating $H_0$ and $\ga$ is given by $r_{\B{G},\Fr{a}}^{-1}(r_{\B{G},\Fr{a}}(K_0)) \cap G_0$.     
\end{pf}

We need three results in order to derive \refT{Borel} from \refT{CSS}. The first is a result of Chahal \cite{Chahal80} which establishes the congruence subgroup property for solvable algebraic groups defined over number fields.

\begin{thm}\label{T:Ch} If $\B{S}$ be a solvable $k$--algebraic group, then every finite index subgroup of $\B{S}(\Cal{O}_k)$ is a congruence subgroup.
\end{thm}

The second is due to Borel---see \cite{Borel91}.

\begin{thm}\label{T:BorelField}
If $\B{G}$ is a $k$--linear algebraic group and $\B{B}$ is a Borel subgroup, then $\B{B}$ is a $k^\prime$--defined closed subgroup for some finite extension $k^\prime/k$.
\end{thm}

The final result is due to Mal'cev.

\begin{thm}\label{T:Malcev}
If $\B{S}$ is a linear, solvable $k$--algebraic group, then every subgroup of $\B{S}(\Cal{O}_k)$ is separable in $\B{S}(\Cal{O}_k)$.
\end{thm}

\begin{pf}[Proof of \refT{Borel}] Let $\B{G}$ be a $k$--algebraic group with Borel subgroup $\B{B}$. We are required to separate each subgroup $S$ of $\B{B}(\Cal{O}_k)$ in $\B{G}(\Cal{O}_k)$. By \refT{BorelField} and \refT{Ch}, $\B{B}$ is a closed $k^\prime$--algebraic subgroup of $\B{G}$ for some finite extension $k^\prime/k$ and $\B{B}(\Cal{O}_{k^\prime})$ has the congruence subgroup property. Hence, \refT{CSS} is applicable and $S$ is separable in $\B{G}(\Cal{O}_{k^\prime})$ if and only if $S$ is separable in $\B{B}(\Cal{O}_{k^\prime})$. By \refT{Malcev}, $S$ is separable in $\B{B}(\Cal{O}_{k^\prime})$ and hence separable in $\B{G}(\Cal{O}_{k^\prime})$. Therefore, by \refL{Basic} (d), $S$ is separable in $\B{G}(\Cal{O}_k)$.
\end{pf}

\section{Corollaries to \refT{Borel}}

Our first corollary shows the conclusions of \refT{Borel} hold for any subgroup of $\B{G}$ commensurable with $\B{G}(\Cal{O}_k)$; we call such subgroups \emph{$k$--arithmetic subgroups}.

\begin{cor}\label{C:CBorel}
Let $\B{G}$ be a $k$--algebraic group, $\La$ a $k$--arithmetic subgroup in $\B{G}$, and $\B{B}$ a Borel subgroup of $\B{G}$. Then every subgroup of $\La \cap \B{B}$ is separable in $\La$.
\end{cor}

\begin{pf} For a subgroup $S$ of $\B{B}\cap \La$, by \refL{Basic} (c), it suffices to separate $S \cap \B{G}(\Cal{O}_k)$ in $\B{G}(\Cal{O}_k) \cap \La$. Since $S \cap \B{G}(\Cal{O}_k)$ is a subgroup of $\B{B}(\Cal{O}_k)$, by \refT{Borel}, $S \cap\B{G}(\Cal{O}_k)$ is separable in $\B{G}(\Cal{O}_k)$. Thus, $S \cap \B{G}(\Cal{O}_k)$ is separable in $\B{G}(\Cal{O}_k) \cap \La$.
\end{pf}  

\begin{cor}\label{C:3Cor}
Let $\B{G}$ be a $k$--algebraic group and $S$ is a solvable subgroup of $\B{G}(\Cal{O}_k)$. Then $S$ is separable in $\B{G}(\Cal{O}_k)$. 
\end{cor}

\begin{pf} Since every solvable subgroup is virtually contained in a Borel subgroup (see \cite[p. 137]{Borel91}), by \refL{Basic} (c), it suffices to separate $S \cap \B{B}$ in $\B{G}(\Cal{O}_k)$. The latter is done using \refT{Borel}.  
\end{pf}

The following is a consequence of \cite{Scott78}.

\begin{thm}\label{T:GImmersion} If $\rho\co N \lra M$ is a $\pi_1$--injective immersion of an infrasolv manifold $N$ into an arithmetic $\B{G}$--orbifold $M$, then there exists a finite cover $\psi\co M^\prime \lra M$ such that $\rho$ lifts to an embedding. 
\end{thm}

\section{Applications to $X$--hyperbolic manifolds}

In this section we prove \refC{Sep} and corollaries specific to lattices in the isometry group of $X$--hyperbolic space. 

\begin{pf}[Proof of \refC{Sep}] \refC{Sep} requires for each $v \in \prt \Hy_Y^n$ and arithmetic lattices $\La$ in $\Isom(\Hy_Y^n)$, we separate each subgroup of 
\[ \tri_v(\La) = \Stab(v) \cap \La \] 
in $\La$, where $Y=\R$, $\C$, $\BB{H}$, or $\BB{O}$. For this, we split our consideration naturally into three cases depending on the stabilizer trichotomy for the groups $\tri_v(\La)$ given at the end of previous chapter. Since $X$--hyperbolic lattices are residually finite it follows easily from \refL{Basic} (c) that subgroups in case (1) are separable. For $X=\R$ or $\C$, case (2) follows exactly the proof in \cite{Hamilton01} on noting $\GL(n;\C) \lra \GL(2n;\R)$. For $X=\BB{H}$ or $\BB{O}$, since every lattice in $\Isom(\Hy_\BB{H}^n)$ and $\Isom(\Hy_\BB{O}^2)$ is arithmetic (see \cite{Corlette92} and \cite{GromovSchoen92}), we can apply \refC{3Cor} to separate. For (3), as peripheral subgroups are virtually nilpotent, \refC{3Cor} handles this case. To be complete, we first realize the arithmetic lattice $\La$ as a subgroup of $\GL(m;\Q)$ with a finite index subgroup in $\GL(m;\Z)$ and finish by applying \refC{CBorel} with \refC{3Cor}. 
\end{pf}

\begin{rem}
\cite{Hamilton01} proved that in a cocompact lattice $\La$ of $\Isom(\Hy_\R^n)$, every virtually abelian subgroup is separable. As her proof does not require arithmeticity, our proof of \refC{Sep} uses arithmeticity only in (3). 
\end{rem}

We conclude this chapter with a geometric corollary of particular interest in the topological classification of cusp cross-sections of arithmetic $X$--orbifolds. 

\begin{thm}\label{T:Cusp}
\begin{description}
\item[(a)] A flat $n$--manifold is diffeomorphic to a cusp cross-section of an arithmetic real hyperbolic $(n+1)$--orbifold if and only if
$\pi_1(M^n)$ injects into an arithmetic real hyperbolic $(n+1)$--lattice.
\item[(b)] An almost flat $(2n+1)$--manifold  $M^{2n+1}$ modelled on $\Fr{N}_{2n+1}$ is diffeomorphic to a cusp cross-section of an arithmetic complex hyperbolic $(n+1)$--orbifold if and only if $\pi_1(M^{2n+1})$ injects into an arithmetic complex hyperbolic $(n+1)$--lattice.
\item[(c)] An almost flat $(4n+3)$--manifold $M^{4n+3}$ modelled on $\Fr{N}_{4n+3}(\BB{H})$ is diffeomorphic to a cusp cross-section of a quaternionic hyperbolic $(n+1)$--orbifold if and only if $\pi_1(M^{2n+1})$ injects into a quaternionic hyperbolic $(n+1)$--lattice. 
\item[(d)] An almost flat $15$--manifold $M^{15}$ modelled on $\Fr{N}_{15}(\BB{O})$ is diffeomorphic to a cusp cross-section of an octonionic hyperbolic $16$--orbifold if and only if $\pi_1(M^{15})$ injects into an octonionic hyperbolic $16$--lattice. 
\end{description}
\end{thm}

\chapter{Cusps of $X$--hyperbolic manifolds}

One goal of this chapter is to give a classification of cusp cross-sections of arithmetic $X$--hyperbolic $n$--orbifolds. By
\refT{Cusp}, we are reduced to classifying AB-groups admitting injections into arithmetic $X$--hyperbolic lattices. The main point is to prove this is equivalent admitting injections into arithmetically defined subgroups of unitary affine groups. The latter groups are easier to work with in regard to this problem, as the generalized Bieberbach theorems ensure the existence of injections. The proof of this reduction relies on being able to realize unitary affine groups as real algebraic subgroups of the isometry group of $X$--hyperbolic space.  In total, this is straightforward with the bulk of the material consisting of terminology, notation, and formal manipulation. We hope the main point is not lost in this.   
 
\section{Algebraic structure of unitary affine groups}

Recall for each maximal compact subgroup $M(X)$ of $\Aut(\Fr{N}_{\ell n-1})$, we defined the unitary affine group $U_M(n-1;X)$ to be $\Fr{N}_{\ell n-1}(X) \rtimes M(X)$. The algebraic structure of these groups is completely determined by the algebraic structure of the maximal compact subgroup. Specifically, $U_M(n-1;X)$ is $k$--algebraic if and only if $M$ is $k$--algebraic. In turn, the algebraic structure of $M$ is controlled by the finite index subgroup $\Uni(H_M)$. For these groups, $\Uni(H_M)$ is $k$--algebraic if and only if $H_M$ is defined
over $k$.

In the real setting, these groups are of the form $\Ort(B_\iny)$, where $B_\iny$ is a symmetric, positive definite bilinear form and the form $B_\iny$ will be defined over a subfield $k$ of $\R$. In the complex setting, these groups are of the form $\Uni(H_\iny)$, where $H_\iny$ is a hermitian form of signature $(n-1,0)$ and $H_\iny$ will be defined over a subfield $k$ of $\C$. In the quaternionic setting, these groups are of the form $\Uni(H_\iny)$, where $H_\iny$ is a positive definite hermitian form and $H_\iny$ will be defined over a subalgebra $A$ of $\BB{H}$. Our only interest is when $k$ is a number field in the first two settings or $A$ is a quaternion algebra defined over a number field in the last setting. 

For an AB-group $\Ga$ modelled on $\Fr{N}_{\ell n-1}(X)$, $\Ga$ is conjugated into a subgroup of a unitary affine group $U_M(n-1;X)$ for any $M(X)$. If this unitary affine group is $k$--algebraic and $\Ga$ is contained in the $k$--points, we say $\Ga$ is \emph{$k$--defined}. When $\Ga$ is commensurable with the $\Cal{O}_k$--points ($\Cal{O}_k$ is either the ring of integers of $k$ or a $\Cal{O}_k$--order in the quaternion algebra), we say $\Ga$ is a \emph{$k$--arithmetic subgroup}. Note if $\Ga$ is $k$--defined, then by conjugating by a Heisenberg dilation, we can arrange for $\Ga$ to be commensurable with a subgroup of the $\Cal{O}_k$--points of the unitary affine group. 

In the quaternionic setting, we can realize $U_M(n-1;\BB{H})$ as $\wh{k}$--algebraic subgroup of $\GL(m;\R)$, where $\wh{k}$ is the
field for which the quaternion algebra $A$ is defined. For a $\Cal{O}_{\wh{k}}$--order $\Cal{O}$ in $A$, if $\Ga$ has a finite index subgroup in the $\Cal{O}$--points of some unitary affine group $U_M(n-1;\BB{H})$, when we realize $U_M(n-1;\BB{H})$ as a $\wh{k}$--algebraic group, $\Ga$ will have a finite index subgroup in the $\Cal{O}_{\wh{k}}$--points of this group. 

In our notation, we will refer to $U_M(n-1;\BB{H})$ as being $A$--defined, subgroups $\Ga$ which are commensurable with
$\Uni(n-1;\Cal{O})$ for some $\Cal{O}_{\wh{k}}$--order $\Cal{O}$ as being $A$--arithmetic, and homomorphisms as being $A$--defined. Since when we realize $\Uni(n-1;\BB{H})$ as a $\wh{k}$--algebraic group, these definitions correspond to the standard algebraic definitions (over the field $\wh{k}$), this is only a slight abuse of notation.

\section{$k$--monomorphisms of unitary affine groups into the isometry group} 

We start by characterizing when a unitary affine group admits a $k$--algebraic structure via embeddings into the isometry group of $X$--hyperbolic space. 

If $U_M(n-1;X)$ is a $k$--algebraic unitary affine group, then $H_{M(X)}$, the associated hermitian form for $M(X)$, is defined over $k$. Set $H = H_{M(X)} \op D_2$, with $H$ defined on $X^{n-1} \op X^2$ and $(X^2,D_2)$ is a $k$--defined $X$--hyperbolic plane. Finally, let $V_\iny$ denote the $H$--orthogonal complement in $X^{n+1}$ of a pair of $X$--linearly independent, $k$--defined, $H$--null vectors $v$ and $v_0$ in $(X^2,D_2)$.

Any isometric $k$--isomorphism 
\[ \psi\co (X^{n-1},H_{M(X)}) \lra (V_\iny, H_{|V_\iny}) \] 
induces a $k$--isomorphism 
\[ \rho\co U_M(n-1;X) \lra MN, \] 
where $N$ and $M$ are factors in the Iwasawa decomposition induced on $\Stab(v)$ with respect to the above pair of $H$--null vectors. Since both vectors are $k$--defined, it follows that $MN$ is $k$--algebraic. Consequently, we have the following proposition.

\begin{prop}\label{P:311}
$U_M(n-1;X)$ is a $k$--algebraic group if and only if there exists a hermitian form $H$ of signature $(n,1)$ defined over $k$ and a $k$--isomorphism 
\[ \rho\co U_M(n-1;X) \lra MN < \Isom(\Hy_X^n) \] 
where $\Hy_X^n$ is modelled on $H$.
\end{prop}

\section{A necessary and sufficient condition for arithmeticity}

We are ready to classify cusp cross-sections of arithmetic hyperbolic lattices. Above, we related the algebraic structure of abstractly defined unitary affine groups via embeddings into the isometry group of $X$--hyperbolic space. We now do the same for AB-groups which is achieved with the following proposition. 

\begin{prop}\label{P:321}
$\Ga$ is a $k$--defined AB-group modelled on $\Fr{N}_{\ell n-1}(X)$ if and only if there exists a $k$--defined hermitian form $H$ modelling
$X$--hyperbolic $n$--space, a subgroup $\La$ of $\Uni(H;k)$ commensurable with $\Uni(H;\Cal{O}_k)$, and an injection 
\[ \rho\co \Ga \lra \Stab(v) \cap \La \] 
for some $k$--defined $H$--null vector $v$.
\end{prop} 

\begin{pf} For the direct implication, assume $\Ga$ is contained in a $k$--defined unitary affine group $U_M(n-1;X)$, and let 
\[ \rho\co U_M(n-1;X) \lra MN < \Uni(H) \] 
be $k$--isomorphism given by \refP{311}. This provides us with a $k$--monomorphism 
\[ \rho\co U_M(n-1;X) \lra \Uni(H) \] 
of $k$--algebraic groups. Hence, by \refL{323}, there exists $\La$ in $\Uni(H;k)$, commensurable with $\Uni(H;\Cal{O}_k)$ such that $\rho(\Ga)$ is contained in $\La$, as asserted. 

For the reverse implication, we assume the existence of $H$, $\La$, $\rho$, and $v$. Note for the Fitting subgroup $L$ of $\Ga$, $\rho(L)$ is contained in some nilpotent factor $N$ of an Iwasawa decomposition. The nilpotent group $N$ is $k$--algebraic, since $L$ is Zariski dense and consists of $k$--points. Since $\rho(\Ga)$ is virtually contained in $N$, $\rho(\Ga)$ is contained in $MN$ for the compact factor $M$ of an Iwasawa decomposition $MAN$ of $\Stab(v)$. Since the group $M$ can be selected to be $k$--algebraic, we have $\rho(\Ga)$ is contained in $MN$, where $MN$ is a $k$--algebraic unitary affine group, as desired.
\end{pf}

For an AB-group $\Ga$ modelled on $\Fr{N}_{\ell n-1}(X)$, we say $\Ga$ is \emph{arithmetically admissible} if there exists an 
arithmetic $X$--hyperbolic $n$--lattice $\La$ such that $\Ga$ is isomorphic to $\tri_v(\La)$. Altogether we have the following theorem which classifies the arithmetically admissible AB-groups (part (a) is proved in \cite{LongReid02}). 

\begin{thm}[Cusp classification theorem]\label{T:324}
Let $\Ga$ be an AB-group modelled on $\Fr{N}_{\ell n-1}(X)$.
\begin{description}
\item[(a)] For $X=\R$, $\Ga$ is arithmetically admissible if and only if $\Ga$ is a $\Q$--arithmetic subgroup in $\R^{n-1} \rtimes \Ort(B_\iny)$, where  $B_\iny$ is a $\Q$--defined, positive definite, symmetric bilinear form  on $\R^{n-1}$. 
\item[(b)] For $X=\C$, $\Ga$ is arithmetically admissible if and only if $\Ga$ is a $k$--arithmetic subgroup in a unitary affine group for some imaginary quadratic number field $k$.
\item[(c)] For $X=\BB{H}$, $\Ga$ is arithmetically admissible if and only if $\Ga$ is a $A$--arithmetic subgroup in a unitary affine group, for some ramified quaternion $\Q$--algebra $A$.   
\end{description}
\end{thm}

\begin{pf} The direct implication is immediate in all three case. For the converse, assume $\Ga$ is a $k$--arithmetic subgroup in a unitary affine group, where $k$ is as above. By \refP{321}, there exists a $k$--defined hermitian form $H$ modelling $X$--hyperbolic $n$--space, a subgroup $\La$ of $\Uni(H;k)$ commensurable with $\Uni(H;\Cal{O}_k)$, and an injection 
\[ \rho\co \Ga \lra \Stab(v) \cap \La \] 
for some $k$--defined light-like vector $v$. $\rho(\Ga)$ must be a finite index subgroup of $\tri_v(\La)$ and by \refT{ClassicationNCL}, $\La$ is an arithmetic subgroup. In this injection we cannot ensure that $\rho(\Ga) = \tri_v(\La)$. As $\La$ is an arithmetic subgroup in the $k$--algebraic group $\Uni(H)$, by \refC{Sep}, we can find a finite index subgroup $\Pi$ of $\La$ such that $\rho(\Ga) = \tri_v(\Pi)$. Specifically, select a complete set of coset representatives for $\tri_v(\La)/\rho(\Ga)$, say $\al_1,\dots,\al_r$. By \refC{Sep}, there exists a finite index subgroup $\Pi$ of $\La$ such that $\rho(\Ga)$ is contained in $\Pi$ and for each $j=1,\dots,r$, $\al_j \notin \Pi$. Therefore $\tri_v(\Pi) = \rho(\Ga)$, since $\Pi \cap \tri_v(\La) = \rho(\Ga)$.
\end{pf}

\section{The holonomy theorem}

For an AB-group $\Ga$ modelled on $\Fr{N}_{2n-1}$, we say the holonomy group $\te$ of $\Ga$ is \emph{complex} if $\te$ is contained in $\Uni(H_{M(X)})$ for the holonomy representation, and otherwise say $\te$ is \emph{anticomplex}. We have the following alternative characterization based on the structure of the holonomy representation.  

\begin{cor}[Holonomy theorem]\label{C:Holonomy}
If $\Ga$ is an AB-group modelled on $\Fr{N}_{2n-1}$ with complex holonomy, then $\Ga$ is arithmetically admissible if and only if the
holonomy representation $\vp$ is conjugate to a representation into $\GL(n-1;k)$ for some imaginary quadratic number field.
\end{cor}

\begin{pf} If $\Ga$ is arithmetically admissible, then from \refT{324}, there exists a $k$--defined unitary affine group $U_M(n-1;k)$ such that $\Ga$ is conjugate into $U_M(n-1;k)$ and commensurable with $U_M(n-1;\Cal{O}_k)$, for some imaginary quadratic number field $k$. This yields an injective homomorphism 
\[ \rho\co \te \lra M(k). \] 
Since $\te$ complex, $\rho(\te)$ is contained in $\Uni(H_{M(X)};k)$, which is a subgroup of $\GL(n-1;k)$, as desired.

For the converse, assume the holonomy representation of $\te$ maps into $\GL(n-1;k)$, for some imaginary quadratic number field $k$. By taking the $\te$--average of any $k$--defined hermitian form, we see this representation is contained in a $k$--defined unitary group $\Uni(H_{M(X)};k)$. Specifically, for any $k$--defined hermitian form $h$, define the $\te$--average of $h$ to be the $k$--defined hermitian form given by
\[ h_\te(z,w) = \sum_{\ga \in \te} h(\ga z,\ga w). \]
Using this representation and a presentation for $\Ga$, we get a system of linear homogenous equations with coefficients in $k$. Since $\rho$ is conjugate to the holonomy representation, by the generalized Bieberbach theorems, this system has a solution which yields a faithful representation into $\Fr{N}_{2n-1}(k) \rtimes \Uni(H_{M(X)};k)$. Conjugating by a Heisenberg dilation to ensure the Fitting subgroup consists of $k$--integral entries, we see $\Ga$ is $k$--arithmetic. Therefore, by \refT{324}, $\Ga$ is arithmetically admissible.
\end{pf}

\section{Density of cusp shapes}

The proof of \refT{SDensity} is established in three steps. We first show each representation $\rho$ in $\Cal{R}(\Ga;\Q)$ produces a similarity class that arises in an arithmetic real hyperbolic $(n+1)$--orbifold. Next, we demonstrate the density of $\Cal{R}(\Ga;\Q)$ in $\Cal{R}_f(\Ga)$. Finally, we produce a continuous surjective map from $\Cal{R}_f(\Ga)$ to both $\Cal{F}_f(\Ga), \Cal{S}_f(\Ga)$. 

Our first proposition follows from \refT{324} and provides the first step in the quest for density.

\begin{prop}\label{P:StrongClass}
Let $\rho \in \Cal{R}(\Ga;\Q)$ and $\Ort_K(n)$ be any $\Q$--defined orthogonal affine group where $\rho(\Ga)$ is a $\Q$--arithmetic subgroup of $\Ort_K(n)$. Then there exists an arithmetic hyperbolic $(n+1)$--lattice $\La$ in $\Isom(\Hy^n)$ and an injection $\psi\co \Ort_K(n) \lra \Isom(\Hy^n)$ such that $\psi(\rho(\Ga))$ is a maximal peripheral subgroup of $\La$.
\end{prop}

\begin{prop}\label{P:Orbit}
For every Bieberbach group $\Ga$, $\Cal{R}(\Ga;\Q)$ is dense in $\Cal{R}_f(\Ga)$.
\end{prop}

To prove \refP{Orbit}, we need a pair of auxiliary lemmas.

\begin{lemma}\label{L:DenseOrbit}
If $G$ is a topological group with dense subgroup $H$ and $X$ is a topological space with a continuous transitive $G$--action, then the $H$--orbit of any $x$ in $X$ is dense in $X$.
\end{lemma}

\begin{pf} Given a dense subgroup $H$ of $G$ and $x$ in $X$, it suffices to show for each $y$ in $X$, there exists a net $h_\al$ in $H$ such that $h_\al \cdot x$ converges to $y$. By assumption $G$ acts transitively on $X$ and so there exists $g$ in $G$ such that $g\cdot x = y$. Since $H$ is dense, there exists a net $h_\al$ in $H$ such that $\lim h_\al = g$. Finally, since the action map $G \times X \lra X$ is continuous, it follows that
\[ \lim h_\al \cdot x = (\lim h_\al) \cdot x = g \cdot x = y. \]
\end{pf}

\begin{lemma}\label{L:QPoints} For every crystallographic group $\Ga$, $\Cal{R}(\Ga;\Q)$ is nonempty.
\end{lemma}

\begin{pf} Let $T_\Ga$ be the maximal translational subgroup and $\te$ the holonomy. We can view $T_\Ga = \Z^n$ and thus the holonomy representation becomes 
\[ \vp\co \te \lra \GL(n;\Z). \] 
After conjugating by a dilation if necessary this yields a faithful representation 
\[ \rho\co \Ga \lra \Z^n \rtimes \GL(n;\Z). \] 
To obtain a $\Q$--defined $\te$--invariant, positive definite bilinear form $B_\te$, we simply take the $\te$--average of any $\Q$--defined positive definite bilinear form $B$. Specifically, define the $\te$--average of $B$ to be
\[ B_\te(x,y) = \frac{1}{\abs{\te}}\sum_{\ga \in \te} B(\ga x,\ga y). \]
It is a simple matter that $B_\te$ is $\Q$--defined, positive definite, and $\te$--invariant and $\Ga$ is a finite index subgroup of $\Z^n \rtimes \Ort(B_\te;\Z)$.
\end{pf}

\begin{pf}[Proof of \refP{Orbit}] As expected, we seek to apply \refL{DenseOrbit}, and must ensure the conditions are satisfied by $X=\Cal{R}_f(\Ga)$ and $G=\Aff(n)$. To begin, the topology on $\Cal{R}_f(\Ga)$ is the subspace topology induced by viewing $\Cal{R}_f(\Ga)$ as a subspace of the $\Aff(n)$--representation space. Visibly, the $\Aff(n)$--action on the representation space is continuous, and so by restriction the action of $\Aff(n)$ on $\Cal{R}_f(\Ga)$ is continuous. Less obvious is the transitivity of the $\Aff(n)$--action on $\Cal{R}_f(\Ga)$. However, this is precisely the statement of one part of the Bieberbach theorems. Thus, by \refL{DenseOrbit}, for $H=\Q^n \rtimes \GL(n;\Q)$ and any $\rho$ in $\Cal{R}_f(\Ga)$, the $H$--orbit of $\rho$ is dense in $\Cal{R}_f(\Ga)$. We assert for each $\al$ in $H$ and $\rho \in \Cal{R}(\Ga;\Q)$, the conjugate representation $\mu_\al \circ \rho$ is in $\Cal{R}(\Ga;\Q)$, where $\mu_\al(\la) = \al^{-1}\la\al$. To see this, let $\Ort_K(n)$ be a $\Q$--defined orthogonal affine group for which $\rho(\Ga)$ is a $\Q$--arithmetic subgroup of $\Ort_K(n)$. Conjugation by $\al$ yields an isomorphism  
\[ \mu_\al\co \textrm{O}_K(n) \lra \textrm{O}_{\be^{-1}K\be}(n) \] 
where $\be$ in $\GL(n;\Q)$ is the linear factor (or second coordinate) for $\al$. Since $\be$ resides in $\GL(n;\Q)$, the symmetric, positive definite form associated to $\be^{-1}K\be$ is $\Q$--defined, being $\Q$--equivalent to the $\Q$--defined form $B_K$. Moreover, this isomorphism between $\Ort_K(n)$ and $\Ort_{\be^{-1}K\be}(n)$ is $\Q$--defined. Therefore, by \refL{323}, any $\Q$--arithmetic subgroup of $\Ort_K(n)$ is mapped to a $\Q$--arithmetic subgroup of $\Ort_{\be^{-1}K\be}(n)$, and thus $\mu_\al \circ \rho(\Ga)$ is a $\Q$--arithmetic subgroup of a $\Q$--defined orthogonal affine group as asserted. Of course, this is only useful if $\Cal{R}(\Ga;\Q)$ is nonempty, and by \refL{QPoints} it is. Consequently, there exists a dense $H$--orbit of representations which, by the argument above, resides in $\Cal{R}(\Ga;\Q)$. 
\end{pf}

\begin{pf}[Proof of \refT{SDensity}] The reduction of $\Cal{R}_f(\Ga)$ by the $\Euc(n)$--conjugate action yields a space containing $\Cal{F}_f(\Ga)$, and there exists a continuous surjective map
\[ \Cal{L}\co \Cal{R}_f(\Ga)/\Euc(n) \lra \Cal{F}_f(\Ga) \]
given as follows. For $\rho$ in $\Cal{R}_f(\Ga)$, as $\rho(\Ga)$ is an $\Aff(n)$--conjugate of a Bieberbach group, by the Bieberbach theorems $\rho(\Ga)$ projects to a finite group $\te$ in $\GL(n;\R)$. Taking the $\te$--average
\[ B_\te(x,y) = \frac{1}{\abs{\te}}\sum_{g \in \te} \innp{gx,gy} \]
of the standard inner product $\innp{\cdot,\cdot}$ on $\R^n$ produces a maximal compact subgroup $K=\Ort(B_\te)$ such that $\rho(\Ga)$ is contained in the orthogonal affine group $\Ort_K(n)$. Up post-composition with an inner automorphism of $\Euc(n)$, there exists a unique $S_\rho$ in $\GL(n;\R)$ conjugating $\Ort_K(n)$ to $\Euc(n)$. From this, we define $\Cal{L}(\rho) = S^{-1}_\rho \rho S_\rho$.

If $\rho_n$ is a sequence of representations in $\Cal{R}_f(\Ga)/\Euc(n)$ converging to $\rho$ in $\Cal{R}_f(\Ga)/\Euc(n)$, the sequence $\Cal{L}(\rho_n)$ converges to  $\Cal{L}(\rho)$ in $\Cal{F}_f(\Ga)$. For a free abelian group this is immediate since $\te$ is trivial. For nontrivial $\te$, this follows from the convergence of the maximal compact subgroups $K$ arising from the $\te$--average---a small change in the image of $\te$ results in a small change in the $\te$--average. We briefly explain this. To begin, the maximal compact subgroup $K$ depends only on $\te$. For a sequence of representations $\rho_n$ converging to $\rho$, let $\te_n$ be the image of $\rho_n(\Ga)$ under projection onto $\GL(n;\R)$. It follows the groups $\te_n$ converge to $\te$ and thus $B_{\te_n}$ converges to $B_\te$ in the space of positive definite, symmetric matrices---we use the standard basis to associate these matrices to the $\te$--average forms. From this we see the maximal compact subgroups $K_n$ converge to $K$. The conjugating matrices $S_n$ need not converge to the conjugating matrix $S$ for $\rho$. However, up to left multiplication in $\Ort(n)$, the sequence does converge and so the sequence of $S_n^{-1}\rho_nS_n$ converges to $S^{-1}\rho S$ in $\Cal{F}_f(\Ga)$. As this is the sequence $\Cal{L}(\rho_n)$, we see $\Cal{L}(\rho_n)$ converges to $\Cal{L}(\rho)$ in $\Cal{F}_n(\Ga)$.   

The desired set of flat similarity classes is the image of $\Cal{L}(\Cal{R}(\Ga;\Q))$ under the projection map
\[ \textrm{Pr}\co \Cal{F}_f(\Ga) \lra \Cal{S}_f(\Ga). \]
That this subset is dense is a consequence of the continuity and surjectivity of $\Cal{L}$ in combination with \refP{Orbit}. For the former, if $\rho(\Ga)$ resides in $\Euc(n)$, the $\te$--average of the standard form is the standard form and thus produces $\Ort(n)$. In particular, we can take $S_\rho = I_n$ and thus $\Cal{L}$ restricted to $\Cal{F}_f(\Ga)$ is the identity. It remains to show each similarity class $\textrm{Pr}(\Cal{L}(\Cal{R}(\Ga;\Q)))$ does arise in a cusp cross-section of an arithmetic hyperbolic $(n+1)$--orbifold. By \refP{StrongClass}, for each $\rho \in \Cal{R}(\Ga;\Q)$ with associated $\Q$--defined orthogonal affine group $\Ort_K(n)$, there exists a faithful representation
\[ \psi\co \Ort_K(n) \lra \Isom(\Hy^{n+1}) \]
and an arithmetic lattice $\La$ such that $\psi(\rho(\Ga))$ is a maximal peripheral subgroup of $\La$. For the flat structure on $\R^n/\Ga$ coming from $\Ort_K(n)$ and the flat structure on the cusp cross-section associated to $\psi(\rho(\Ga))$, this produces a similarity of this pair of flat manifolds. To obtain this for the associated class in $\textrm{Pr}(\Cal{L}(\rho))$, we argue as follows. It could be that the $\Q$--form $B_K$ for $K$ is not the $\te$--average of $\rho(\Ga)$. If this is the case, simply replace $K$ by $\Ort(B_\te)$, and notice this too is a $\Q$--defined orthogonal affine group for which $\rho(\Ga)$ is a $\Q$--arithmetic subgroup. Let $M^\prime$ be the associated flat manifold with this similarity class associated to $\rho$ viewed as a representation into $\Ort_{\Ort(B_\te)}(n)$. Making the same argument as before, we see $M^\prime$ occurs as a cusp cross-section of an arithmetic hyperbolic $(n+1)$--orbifold. By construction, the flat manifold $M^{\prime\prime}$ with similarity class $\textrm{Pr}(\Cal{L}(\rho))$ is similar to $M^\prime$. Hence, every class in the dense subset $\textrm{Pr}(\Cal{L}(\Cal{R}(\Ga;\Q)))$ arises in a cusp cross-section of an arithmetic real hyperbolic $(n+1)$--orbifold. 
\end{pf}

For almost flat manifolds modelled on the $\Fr{N}_{2n-1}$ or $\Fr{N}_{4n-1}(\BB{H})$ an identical argument can be made. However, the associated set $\Cal{R}(\Ga;\Q)$ need not be nonempty, and consequently there is a dichotomy. 

\begin{thm}\label{T:AlmostFlat}
\begin{description}
\item[(a)] For an almost flat $(2n-1)$--manifold $N$ modelled on $\Fr{N}_{2n-1}$, the realizable almost flat similarity classes in the cusp cross-sections of arithmetic complex hyperbolic $n$--orbifolds is either empty or dense in the space of almost flat similarity classes.
\item[(b)] For an almost flat $(4n-1)$--manifold $N$ modelled on $\Fr{N}_{4n-1}(\BB{H})$, the realizable almost flat similarity classes in the cusp cross-sections of quaternionic hyperbolic $n$--orbifolds is either empty or dense in the space of almost flat similarity classes.
\end{description}
\end{thm}

Together with \refT{Nil3Thm}, we obtain:

\begin{cor} For a $\Nil$ 3--manifold $N$, the $\Nil$ similarity classes that arise in the cusp cross-sections of arithmetic complex hyperbolic $2$--orbifolds are dense in the space of $\Nil$ similarity classes.
\end{cor}

This also holds for $\Sol$ 3--manifolds---see the next chapter for more on this.

\begin{cor}\label{C:SolDensity} For a $\Sol$ 3--manifold $S$, the $\Sol$ similarity classes that arise in the cusp cross-sections of generalized Hilbert modular surfaces are dense in the space of $\Sol$ similarity classes.
\end{cor}

Finally, with the easily established converse of \refP{StrongClass}, we obtain the following geometric classification theorem for cusp cross-sections of arithmetic hyperbolic $(n+1)$--orbifolds.

\begin{thm} For a flat $n$--manifolds, $\textrm{Pr}(\Cal{L}(\Cal{R}(\pi_1(M);\Q)))$ is precisely the set of flat similarity classes on $M$ that arise in cusp cross-sections of arithmetic hyperbolic $(n+1)$--orbifolds.
\end{thm}

This persists in the complex and quaternionic hyperbolic settings.

\section{Orbifold to manifold promotion}

In this short section we prove the following result.

\begin{thm}\label{T:Density}
For the $n$--torus, the realizable flat similarity classes in the cusp cross-sections of arithmetic real hyperbolic $(n+1)$--manifolds are dense in the space of flat similarity classes.
\end{thm}

\refT{Density} is a consequence of the following proposition whose proof is essentially a reproduction of Borel's proof of Selberg's lemma \cite{Borel63}.

\begin{prop}\label{P:UnipotentSelberg} If $k/\Q$ is a finite extension and $\La$ a finitely generated subgroup of $\GL(n;k)$ with unipotent subgroup $\Ga$, then there exists a torsion free, finite index subgroup $\La_0$ of $\La$ such that $\Ga$ is contained in $\La_0$.
\end{prop}

\begin{pf} Let $\la_1,\dots,\la_r$ be a finite generating set for $\La$, $c_{i,j,\ell}$ be the $(i,j)$--coefficient of $\la_\ell$, and $R$ be the subring of $k$ generated by $\set{c_{i,j,\ell}}$. By assumption, $\Ga$ is conjugate in $\GL(n;\C)$ into the group of upper triangular matrices with ones along the diagonal. In particular, the characteristic polynomial $p_\ga(t)$ for each $\ga$ in $\Ga$ is $(t-1)^n$. For any torsion element $\eta$ in $\La$, the characteristic polynomial $p_\eta(t)$ of $\eta$ has only roots of unity for its zeroes. Since $k/\Q$ is a finite extension and $n$ is fixed, there are only finitely many degree $n$ monic polynomials in $k[t]$ having only roots of unity for their roots. Let $p_1(t),\dots,p_s(t)$ denote those monic polynomials with coefficients in $R$ with this property. As our concern is solely with nontrivial torsion elements, we further insist each of the polynomials has a root distinct from $1$. For each such polynomial $p_j(t)$, there are finitely many prime ideals $\Fr{p}$ of $R$ such that $(t-1)^n = p_j(t)$ modulo $\Fr{p}$. To see this, we first exclude all prime ideals $\Fr{p}$ in $R$ such that $\ch(R/\Fr{p})\leq n$. Since for each prime $p$ of $\Z$ there are finitely many prime ideals $\Fr{p}$ of $R$ such that $\ch(R/\Fr{p})=p$, this is a finite set. Next, as 
\[ p_j(t) - (t-1)^n = \sum_{m=1}^n \la_{m,j}t^m \] 
is nonzero, there exists $i$ such that $\la_{i,j}$ is nonzero. There are only finitely many prime ideals $\Fr{p}_{j,1},\dots,\Fr{p}_{j,\ell_j}$ such that $\la_{i,j} = 0 \mod \Fr{p}_{j,\ell_k}$, and so for any other prime ideal $\Fr{q}$, it follows $p_j(t) \ne (t-1)^n$ modulo $\Fr{q}$. Excluding this finite collection $\Cal{P}_j$ of prime ideals of $R$, for any selection $\Fr{q} \notin \Cal{P}_j$, we have $p_j(t)$ not equal to $(t-1)^n$ modulo $\Fr{q}$. Repeating this argument for each $j$, we obtain the desired ideal set $\Cal{P}$. For $\Fr{q} \notin \Cal{P}$, no torsion element $r_\Fr{q}(\eta)$ cannot reside in $r_\Fr{q}(\Ga)$ as every element of $r_\Fr{q}(\Ga)$ has characteristic polynomial $(t-1)^n$ and $\eta$ does not share this trait. Therefore, $r_\Fr{q}^{-1}(r_\Fr{q}(\Ga))$ is a torsion free finite index subgroup of $\La$ containing $\Ga$, as sought.
\end{pf}

\begin{pf}[Proof of \refT{Density}] It suffices to show for each $\rho$ in $\Cal{R}(\Z^n;\Q)$ the induced representation given by \refP{StrongClass} is such that $\Z^n$ is contained in a torsion free finite index subgroup of the target lattice $\La$. By construction, the representation $\rho\co\Z^n \lra \La$ maps $\Z^n$ into a unipotent subgroup of $\La$ since the groups $N$ in the Iwasawa decomposition are unipotent. As the target lattice $\La$ is arithmetic, $\La$ is finitely presentable (\cite[Cor. 13.25]{Raghunathan72}) and conjugate into the $k$--points of $\Isom(\Hy^{n+1})$ for some number field $k$. Thus \refP{UnipotentSelberg} is applicable and yields a torsion free finite index subgroup $\La_0$ of $\La$ such that $\Z^n$ is contained in $\La_0$. Note if $\Z^n$ is a maximal peripheral subgroup of $\La$, $\Z^n$ is a maximal peripheral subgroup of $\La_0$. In particular, we can realize the associated flat similarity class $[g]$ for $\rho$ in a cusp cross-section of the associated arithmetic manifold for $\La_0$.
\end{pf}

For those infranil manifold groups realizable as lattices in their associated nilpotent Lie group, we say the associated infranil manifold is a \emph{niltorus}. For niltori modelled on either $\Fr{N}_{2n-1}$ or $\Fr{N}_{4n-1}(\BB{H})$, orbifold density is promoted to manifold density with an identical argument using \refP{UnipotentSelberg}.

\chapter{Hilbert and Hilbert--Blumenthal cusps}

Let $k$ be a totally real number field with $[k:\Q]=n$, $\Cal{O}_k$ the ring of integers of $k$, and $\si_1,\dots,\si_n$ the $n$ real embeddings of $k$. The group $\PSL(2;\Cal{O}_k)$ is an arithmetic subgroup of the $n$--fold product $(\PSL(2;\R))^n$ via the embedding
\[ \xi \lmto (\si_1(\xi),\dots,\si_n(\xi)) \] 
for $\xi \in \PSL(2;\Cal{O}_k)$. The group $\PSL(2;\Cal{O}_k)$ is called \emph{the Hilbert modular group}, and through this embedding, $\PSL(2;\Cal{O}_k)$ acts with finite volume on the $n$--fold product of real hyperbolic planes $(\Hy_\R^2)^n$. More generally, we call any subgroup $\La$ of $\PSL(2;k)$ commensurable with $\PSL(2;\Cal{O}_k)$ a \emph{Hilbert modular group} and the quotients $(\Hy_\R^2)^n/\La$, \emph{Hilbert modular varieties}. When $k$ is a real quadratic number field, these quotients are called \emph{Hilbert modular surfaces}---for more on Hilbert modular surfaces, see \cite{Hirzebruch73} or \cite{vanderGeer88}. 

\section{Cusps of Hilbert modular varieties}

For the product geometry $(\Hy_\R^2)^r$, the Iwasawa decomposition of $(\PSL(2;\R))^r$ is given by taking $r$ independent Iwasawa decompositions in each $\PSL(2;\R)$. The standard decomposition in $\PSL(2;\R)$ is given with
\begin{align*}
\B{A} &= \set{\begin{pmatrix} \al & 0 \\ 0 & \al^{-1} \end{pmatrix} ~:~ \al \in \R^+} \\
\B{N} &= \set{\begin{pmatrix} 1 & \be \\ 0 & 1 \end{pmatrix} ~:~ \be \in \R} \\
\B{K} &= \set{\begin{pmatrix} \cos \te & \sin \te \\ -\sin \te & \cos \te \end{pmatrix} ~:~ \te \in [0,2\pi]}.
\end{align*}
The standard Borel subgroup $\B{B}$ of $\PSL(2;\R)$ is $\B{N} \rtimes \B{A}$, and a Borel subgroup of $(\PSL(2;\R))^r$ is the $r$--fold product of Borel subgroups $\B{B}_1,\dots,\B{B}_r$ of $\PSL(2;\R)$.

Cusps, horospheres, and cusp cross-sections are defined as in the $X$--hyperbolic setting via Iwasawa decompositions of $(\PSL(2;\R))^r$. For the Hilbert modular group $\PSL(2;\Cal{O}_k)$ over a totally real number field $k$, the stabilizer of the boundary point corresponding to the Iwasawa decomposition given by the $r$--fold product of the groups $\B{A},\B{N},\B{K}$ is the peripheral subgroup
\[ \tri = \set{\begin{pmatrix} \al & \be \\ 0 & \al^{-1} \end{pmatrix}~:~\be \in \Cal{O}_k,~\al \in \Cal{O}_{k,+}^\times}. \]
Every peripheral subgroup of $\PSL(2;\Cal{O}_k)$ is conjugate in $\PSL(2;k)$ to a group commensurable with $\tri$. For higher real rank geometries one typically has metric pinching in finite volume, noncompact quotients arising from other proper parabolic subgroups. When these parabolic groups are not Borel subgroups, the associated peripheral subgroups are not stabilizers of points on the boundary but instead higher dimensional simplicial complexes. The metric pinching in Hilbert modular varieties arises only from Borel subgroups, a consequence of the fact our geometry has $\Q$--rank one. For completeness, we establish this claim.

\begin{lemma}
If $\B{P}$ is $\Q$--defined proper parabolic subgroup of $\SL(2;\R)^r$ which intersects $\PSL(2;\Cal{O}_k)$ nontrivially, then $\B{P} \cap \PSL(2;\Cal{O}_k)$ is contained in a $\Q$--defined Borel subgroup.
\end{lemma}

\begin{pf} After applying an element of the symmetric group $\textrm{Sym}(r)$, we can assume the $\Q$--defined parabolic $\B{P}$ of $\B{G}=(\PSL(2;\R))^r$ is of the form
\[ \prod_{j=1}^t \B{B}_j \times \prod_{j=t+1}^r \PSL(2;\R). \]
Conjugating in $\B{G}$, we can further assume each $\B{B}_j$ is a standard Borel subgroup $\B{B}_s$ of $\PSL(2;\R)$. If $\ga \in \B{P} \cap \PSL(2;\Cal{O}_k)$, as $\B{P}$ is proper, then under some embedding $\si_\ell$ of $k$, $\si_\ell(\ga) \in \B{B}_s$. In particular,
\[ \si_\ell(\ga) = \begin{pmatrix} \al & \be \\ 0 & \al^{-1} \end{pmatrix}. \]
For any other embedding $\si_j$ of $k$
\[ \si_j(\ga) = \begin{pmatrix} a & b \\ c & d \end{pmatrix} \]
and for some $\te \in \Gal(k_{gal}/\Q)$, 
\[ \si_\ell(\ga) = \begin{pmatrix} \te(a) & \te(b) \\ \te(c) & \te(d) \end{pmatrix}. \]
Therefore, $c =0$ and so for every embedding $\si_j(\ga) \in \B{B}_s$ as asserted.
\end{pf} 

\section{Classifying cusp cross-sections}

We now turn our eyes to Hilbert modular varieties and the classification of their cusp cross-sections. This is achieved with our next theorem analogous to \refT{324}.

\begin{thm}[Correspondence theorem]\label{T:CT}
If $N$ is a $k$--arithmetic torus bundle, then there exists a faithful representation 
\[ \psi\co \pi_1(N) \lra \tri(\PSL(2;\Cal{O}_k)) \]
such that $\psi(\pi_1(N))$ is a finite index subgroup of $\tri(\PSL(2;\Cal{O}_k))$. Moreover, there exists a finite index subgroup $\La$ of $\PSL(2;\Cal{O}_k)$ such that $\tri(\La) = \psi(\pi_1(N))$.
\end{thm}

We defer the proof of \refT{CT} for the moment in order to prove \refT{Torus}. 

\begin{pf}[Proof of \refT{Torus}] Our task is to verify an $(n,n-1)$--torus bundle $N$ is diffeomorphic to a cusp cross-section of a Hilbert modular variety over $k$ if and only if $\pi_1(N)$ is $k$--arithmetic. For the direct implication, since $N$ is diffeomorphic to a cusp cross-section of a Hilbert modular variety, there exists a Hilbert modular group $\La$ and an isomorphism 
\[ \psi\co \pi_1(N) \lra \tri(\La). \] 
To obtain an injective homomorphism 
\[ \rho\co \pi_1(N) \lra k \rtimes k_+^\times \] 
such that $\rho(\pi_1(N))$ is commensurable with $\Cal{O}_k \rtimes \Cal{O}_{k,+}^\times$, we argue as follows. By conjugating by an element $\ga$ of $\PSL(2;k)$, we can assume that
\[ \ga^{-1}\psi(\pi_1(N))\ga \su \B{B}_s(k) = \set{\begin{pmatrix} \be^{-1} &
\al \\ 0 & \be \end{pmatrix} ~:~ \al \in k,~\be \in k_+^\times}. \]
As $\ga$ is in $\PSL(2;k)$, $\ga^{-1}\La\ga$ remains a Hilbert modular group, and moreover, $\ga^{-1}\psi(\pi_1(N))\ga$ is commensurable with 
\[ \tri(\PSL(2;\Cal{O}_k)) = \set{\begin{pmatrix} \be^{-1} & \al \\ 0 &\be \end{pmatrix}~:~\al\in\Cal{O}_k,~\be\in \Cal{O}_{k,+}^\times}. \] 
To obtain the faithful representation $\rho$, we simply compose $\mu_\ga \circ \psi$ with the isomorphism 
\[ \iota\co \B{B}(k) \lra k \rtimes k_+^\times \] given by 
\[ \iota\pr{\begin{pmatrix} \be^{-1} & \al \\ 0 & \be \end{pmatrix}} = (\al,\be). \]   

For the reverse implication, we apply \refT{CT} and \refT{Rig}. Specifically, let $\La$ be the Hilbert modular group
guaranteed by \refT{CT} and let $N^\prime$ denote an embedded cusp cross-section associated with $\tri(\La)$. As a smooth manifold, $N^\prime$ is of the form $\R^{2n-1}/\tri(\La)$ and by \refT{CT}, we have an isomorphism 
\[ \psi\co \pi_1(N) \lra \pi_1(N^\prime). \] 
Applying Mostow's solvable rigidity \refT{Rig}, we obtain the desired diffeomorphism between $N$ and $N^\prime$.
\end{pf}

In the proof of \refT{CT}, the following lemma is required.

\begin{lemma}\label{L:1Hilbert}
If $N$ is a $k$--arithmetic torus bundle, then there exists an injective homomorphism 
\[ \rho\colon \pi_1(N) \lra \Cal{O}_k \rtimes \Cal{O}_{k,+}^\times. \] 
Moreover, $\rho(\pi_1(N))$ is a finite index subgroups of $\Cal{O}_k \rtimes \Cal{O}_{k,+}^\times$.  
\end{lemma}

\begin{pf} Since $N$ is $k$--arithmetic, we have a faithful representation 
\[ \te\colon \pi_1(N) \lra k \rtimes k_+^\times \]
such that $\te(\pi_1(N))$ is commensurable with $\Cal{O}_k \rtimes \Cal{O}_{k,+}^\times$. Hence, given $(\al,\be)$ in $\te(\pi_1(N))$, we have
for some $m \in \N$, 
\[ (\al + \be\al + \be^2\al + \dots+ \be^{m-1}\al,\be^m) \in \Cal{O}_k \rtimes \Cal{O}_{k,+}^\times. \] 
Consequently, $\be^m$ is in $\Cal{O}_{k,+}^\times$ and thus $\be$ is in $\Cal{O}_{k,+}^\times$. Even so, it may be the case that $(\al,\be)$ is not contained in $\Cal{O}_k \rtimes \Cal{O}_{k,+}^\times$. This is rectified as follows. Select a generating set for $\pi_1(N)$, say
$g_1,\dots,g_u$. For each generator, we have 
\[ \te(g_j) = (\al_j,\be_j) \] 
with $\al_j\in k$ and $\be_j \in \Cal{O}_{k,+}^\times$. Since $k$ is the field of fractions of $\Cal{O}_k$, we can select $\la_j \in \Cal{O}_k$ such that 
\[ (0,\la_j)\te(g_j)(0,\la_j)^{-1} \in \Cal{O}_k \rtimes \Cal{O}_{k,+}^\times. \] 
Note 
\[ (0,\la_j)\te(g_j)(0,\la_j)^{-1} = (\la_j\al_j,\be_j), \] 
and so the second coordinate $\be_j$ is unchanged. Finally, for $\la = \la_1\dots\la_u$, define 
\[ \rho = \mu_{(0,\la)} \circ \te, \] 
where $\mu_{(0,\la)}$ denotes the inner automorphism determined by $(0,\la)$. By construction, $\rho$ is a faithful representation of $\pi_1(N)$ onto a finite index subgroup of $\Cal{O}_k \rtimes \Cal{O}_{k,+}^\times$.
\end{pf}

With \refL{1Hilbert} in hand, we prove \refT{CT}.

\begin{pf}[Proof of \refT{CT}] By \refL{1Hilbert}, we have an injective homomorphism 
\[ \rho\co \pi_1(N) \lra \Cal{O}_k \rtimes \Cal{O}_{k,+}^\times \] 
such that $\rho(\pi_1(N))$ is a finite index subgroup. To obtain the injective homomorphism $\psi$, we compose $\rho$ with the isomorphism 
\[ \iota^{-1}\co \Cal{O}_k \rtimes \Cal{O}_{k,+}^\times \lra \tri(\PSL(2;\Cal{O}_k)) \] 
where 
\[ \iota^{-1}(\al,\be) = \begin{pmatrix} \be^{-1} & \al \\ 0 & \be \end{pmatrix}. \] 
That $\psi$ is faithful and $\psi(\pi_1(N))$ is a finite index subgroup of $\tri(\PSL(2;\Cal{O}_k))$ follow immediately from the properties of $\rho$ and $\iota$. To find the desired subgroup $\La$, we apply the Borel subgroup separability \refT{Borel}.
\end{pf}

\section{A question of Hirzebruch}\label{SS:Hirzebruch}

Let $k$ be a totally real number field, $M<k$ an additive group of rank $n$ (the degree of $k$ over $\Q$), and  $V<\Cal{O}_{k,+}^\times$ a finite index subgroup such that for all $\la \in V$, $\la M \su M$. For each pair $(M,V)$, we define the peripheral group
\[ \tri(M,V) = \set{\begin{pmatrix} \be^{-1} & \al \\ 0 & \be
  \end{pmatrix} ~:~ \al \in M,~\be \in V} < \PSL(2;k). \]
For any Hilbert modular variety $X$, the peripheral groups $\tri(\La)$ of $\pi_1(X)$ are conjugate (in $\PSL(2;k)$) to groups of the form $\tri(M,V)$. In
\cite[p. 203]{Hirzebruch73}, Hirzebruch mentions that it is apparently unknown whether or not every $\tri(M,V)$ can occur as a maximal peripheral subgroup of a Hilbert modular group. The following corollary gives an affirmative answer.

\begin{cor}\label{C:Hirzebruch}
For every pair $(M,V)$, there exists a Hilbert modular group $\La$ such that $\tri(\La) = \tri(M,V)$.
\end{cor}

\begin{pf} As in the proof of \refL{1Hilbert}, we can conjugate $\tri(M,V)$ by an element of the form 
\[ \ga = \begin{pmatrix} \la^{-1} & 0 \\ 0 & \la \end{pmatrix}, \] 
with $\la$ in $\Cal{O}_k$, such that $\ga^{-1} \tri(M,V) \ga$ is contained in $\PSL(2;\Cal{O}_k)$. Since $M$ and $V$ are finite index subgroups of $\Cal{O}_k$ and $\Cal{O}_{k,+}^\times$, respectively, $\ga^{-1}\tri(M,V)\ga$ is a finite index subgroup of $\tri(\PSL(2;\Cal{O}_k))$. Thus there exists a finite index subgroup $\La_1$ of $\PSL(2;\Cal{O}_k)$ such that 
\[ \tri(\La_1) = \ga^{-1}\tri(M,V)\ga. \] 
Hence, for $\La = \ga \La_1 \ga^{-1}$, we have $\tri(\La) = \tri(M,V)$. As $\ga \in \PSL(2;k)$, $\La$ is a Hilbert modular group, as required.
\end{pf} 

\section{A simple criterion for arithmeticity}

In this section, we give a simple criterion for the arithmeticity of $(n,m)$--torus bundles, and so produce the analog of the holonomy \refC{Holonomy}. The need for such a result is practical, as it allows one to establish the arithmeticity of a torus bundle computationally.
 
For an (orientable) $(n,n-1)$--torus bundle $M$, since both the base and fiber are aspherical, we have the short exact sequence induced by the long exact sequence of the fiber bundle
\[ 1 \lra \Z^n \lra \pi_1(M) \lra \Z^{n-1} \lra 1. \]
The action of $\Z^{n-1}$ on $\Z^n$ induces a homomorphism 
\[ \vp\colon \Z^{n-1} \lra \SL(n;\Z) \] 
called the \emph{holonomy representation}. Since peripheral subgroups in Hilbert modular groups have faithful holonomy representation, we assume throughout that $\vp$ is faithful. In particular, we obtain a faithful representation of $\pi_1(M)$ into $\Z^n \rtimes \SL(n;\Z)$.

Of primary importance for us here is the holonomy representation together with any finite presentation yields a homogenous linear system of equations with coefficients in $\Z$. This system arises as follows. For ease, select a presentation of the form 
\[ \innp{x_1,\dots,x_n,\ol{y_1},\dots,\ol{y_{n-1}}~:~R} \]
where $x_1,\dots,x_n$ generate $\Z^m$, $\ol{y_1},\dots,\ol{y_{n-1}}$ are lifts of a generating set $y_1,\dots,y_{n-1}$ for $\Z^{n-1}$, and $R$ is a finite set of relations of the form
\[ x_j\ol{y_k} = \ol{y_k}w_{j,k}, \quad w_{j,k} \in
\innp{x_1,\dots,x_n}. \] 
Using the holonomy representation, we can write
\[ x_j = (a_j,I), \quad \ol{y_j} = (b_j,\vp(y_j)) \in \R^n \rtimes
\SL(n;\R). \] 
Each relation in the presentation yields a linear homogenous equation in the vector variables $a_j$ and $b_j$ (see below for an explicit
example of how these equations arise). Namely, we insert the above forms for $x_j$ and $\ol{y_k}$ into the relation and consider only the
first coordinate. The equations we obtain are of the form
\[ a_j + b_k - \vp(y_k) - v_{j,k} = 0 \]
where $w_{j,k} = (v_{j,k},I)$. That this system has integral solutions which yield faithful representations follows from the fact that $\vp$
is faithful and induces a faithful representation of $\pi_1(M)$ into $\Z^n \rtimes \SL(n;\Z)$.  

The main result of this section is a simple criterion for arithmeticity based on the structure of the holonomy representation. In the statement of the proof, $[k:\Q]=n$ and $\rank \Cal{O}_k^\times = n-1$.

\begin{thm}\label{T:Class} If $M$ is an orientable $(n,n-1)$--torus bundle, then $M$ is diffeomorphic to a cusp cross-section of a Hilbert modular variety defined over $k$ if and only if $\vp = \res_{k/\Q}(\chi)$, for some faithful character $\chi\colon \Z^{n-1} \lra \Cal{O}_{k,+}^\times$, 
where $\vp$ is some holonomy representation. 
\end{thm}

\begin{pf} For the direct implication, since $M$ is diffeomorphic to a cusp cross-section of a Hilbert modular variety, by \refT{Torus}, we have a faithful representation 
\[ \rho\co \pi_1(M) \lra \Cal{O}_k \rtimes \Cal{O}_{k,+}^\times. \] 
By restricting scalars from $k$ to $\Q$, we obtain a faithful representation
\[ \res_{k/\Q}(\rho)\co \pi_1(M) \lra \Z^n \rtimes \SL(n;\Z). \]
The proof is completed by noting the holonomy map induced by this representation is simply $\res_{k/\Q}(\chi)$, where $\chi\co \Z^{n-1} \lra \Cal{O}_{k,+}^\times$ is the holonomy representation induced by the representation $\rho$.

For the converse, we seek a faithful representation 
\[ \rho\colon \pi_1(M) \lra \Cal{O}_k \rtimes \Cal{O}_{k,+}^\times. \] 
Since $[k:\Q]=n$ and $\rank \Cal{O}_k^\times = n-1$, the image of $\pi_1(M)$ is necessarily a finite index subgroup. By assumption, we have a faithful character $\chi\colon \Z^{n-1} \lra \Cal{O}_{k,+}^\times$. We extend this to a faithful representation of $\pi_1(M)$ into $\Cal{O}_k \rtimes \Cal{O}_{k,+}^\times$ as follows. Select a presentation as above for $\pi_1(M)$ with generators $x_1,\dots,x_n$ and $\ol{y_1},\dots,\ol{y_{n-1}}$. Write 
\begin{equation}\label{E:Hom}
x_i = (\al_i,1),~\ol{y_i} = (\ga_i,\chi(y_i)) \in k \rtimes \Cal{O}_{k,+}^\times
\end{equation}
where $\al_i$ and $\ga_i$ are to be determined. Using our presentation for $\pi_1(M)$, we obtain a system of linear homogenous equations $\Cal{L}$ with coefficients in $\Cal{O}_k$. As above, solutions to $\Cal{L}$ yield representations of $\pi_1(M)$ into $k\rtimes \Cal{O}_{k,+}^\times$. We assert there is a solution which yields a faithful representation. To see this, by restricting scalars from $k$ to $\Q$, we obtain a linear system $\res_{k/\Q}(\Cal{L})$ with coefficients in $\Z$. Solutions to the system $\res_{k/\Q}(\Cal{L})$ yield representations of $\pi_1(M)$ into $\Z^n \rtimes \SL(n;\Z)$. Moreover, a solution to $\res_{k/\Q}(\Cal{L})$ which yields a faithful representation is equivalent to a solution of $\Cal{L}$ which yields a faithful representation into $\Cal{O}_k \rtimes \Cal{O}_{k,+}^\times$. That such a solution exists with integral coefficients for $\res_{k/\Q}(\Cal{L})$ follows from the faithfulness of $\res_{k/\Q}(\chi)$ and our discussion in the previous subsection. This yields a solution for $\Cal{L}$ with coefficients in $\Cal{O}_k$ which yields a faithful representation. Therefore, $M$ is $k$--arithmetic, since there exists a faithful representation 
\[ \psi\co \pi_1(M) \lra \Cal{O}_k \rtimes \Cal{O}_{k,+}^\times \] 
such that $\psi(\pi_1(M))$ is a finite index subgroup of $\Cal{O}_k \rtimes \Cal{O}_k^\times$. 
\end{pf}

If the character $\chi$ only maps into $\Cal{O}_k^\times$, the above proof yields a faithful representation 
\[ \rho\colon \pi_1(M) \lra \Cal{O}_k \rtimes \Cal{O}_k^\times. \]

\section{Hilbert--Blumenthal modular varieties}

As it requires no more work, we mention the case when $k$ is any algebraic number field. Our interest is in those subgroups of $\PSL(2;k)$
commensurable with $\PSL(2;\Cal{O}_k)$. As in the case when $k$ is totally real, these groups are arithmetic lattices in $(\PSL(2;\R))^{r_1} \times (\PSL(2;\C))^{r_2}$, where $r_1$ is the number of distinct real embeddings of $k$ and $r_2$ is the number of distinct complex embeddings of $k$. These groups are called \emph{Hilbert--Blumenthal modular groups} and the finite volume quotient manifolds of these groups acting on $(\Hy_\R^2)^{r_1} \times (\Hy_\R^3)^{r_2}$ are called \emph{Hilbert--Blumenthal modular varieties}.  

The cusps of these quotients are $(n,m)$--torus bundles, where $n=[k:\Q]$ and $m = \rank \Cal{O}_k^\times$. For $\si_1,\dots,\si_{r_1}$ and $\tau_1,\dots,\tau_{r_2}$, define 
\[ \wt{\si_j}\co k^\times \lra \R^\times/\innp{\pm 1} \] 
and 
\[ \wt{\tau_j}\co k^\times \lra \C^\times/\innp{\pm 1}. \] 
The product of these maps yields
\[ \rho\co k^\times \lra \prod_{j=1}^{r_1} \R^\times/\innp{\pm 1} \times \prod_{j=1}^{r_2} \C^\times/\innp{\pm 1}. \]
We say $V \su k^\times$ is \emph{positive} if $\rho_{|V}$ is an injective homomorphism. Note when $k$ is totally real such a $V$ consists of totally positive numbers. We say an $(n,m)$--torus bundle is \emph{$k$--defined} if there exists a positive subgroup $V$ in $k^\times$ and a faithful representation of $\pi_1(N)$ into $k \rtimes V$. If in addition, $\pi_1(N)$ is commensurable with $\Cal{O}_k \rtimes \Cal{O}_k^\times$ under this representation, we say $N$ is \emph{$k$--arithmetic}.  

\begin{thm}\label{T:GTorus} A virtual $(n,m)$--torus bundle $N$ is diffeomorphic to a cusp cross-section of a Hilbert-Blumenthal modular variety defined over $k$ if and only if $N$ is $k$--arithmetic.
\end{thm}

The proof of \refT{GTorus} is identical to \refT{Torus}---as is \refT{Class} in this more general setting. As before, for any pair $(n,m)$ with $n >2$ and $m>0$, there exist $(n,m)$--torus bundles which are not diffeomorphic to a cusp cross-section of any Hilbert-Blumenthal modular variety.

We can classify cusp cross-sections of irreducible orbifolds modelled on the $n$--fold product ($n>1$) $\prod_{j=1}^n \Hy_X^{m_j}$ whose associated isometry group is, up to finite index, 
\[ \B{G}_{m_1,\dots,m_n} = \prod_{j=1}^n \Isom(\Hy_X^{m_j}), \]
where $X=\R$, $\C$, or $\B{H}$. Irreducible lattices in $\B{G}_{m_1,\dots,m_n}$ exist if and only if each $m_j=2$ or $3$ and $X=\R$, or $m_j=m_k$ for all $j$ and $k$, and by Margulis' arithmeticity theorem \cite{Margulis91}, these lattices are always arithmetic. 

\section{An obstruction to geometric bounding}

Let $W$ be a 1--cusped Hilbert modular surface $W$ with torsion free fundamental group---we call $W$ a \emph{Hilbert modular manifold} in this case. Similar to the thick-thin decomposition of a real hyperbolic $n$--manifold, $W$ has a decomposition comprised of a compact manifold $\wt{W}$ with boundary $S$ and cusp end $S \times \R^+$. Following Schwartz \cite{Schwartz95} (see also \cite{FarbSchwartz96}), we call the manifold $\wt{W}$ the associated \emph{neutered manifold}, and note $\wt{W}$ is a compact 4--manifold with $\Sol$ 3--manifold boundary. The goal of this section is the establishment of a nontrivial obstruction for this geometric situation. The obstruction is obtained by mimicking the argument of Long--Reid \cite{LongReid00} for flat 3--manifolds. This in combination with a calculation of Hirzebruch bears \refT{GB} from the introduction.

In \cite{Hirzebruch73}, Hirzebruch extended his signature formula to Hilbert modular surfaces. The formula relates the signature of the neutered manifold $\wt{W}$ to a Hirzebruch $L$--polynomial evaluated on the Pontrjagin classes of $\wt{W}$ but with a correction term associated to $\prt \wt{W}$. When $\pi_1(W)$ contains torsion, the elliptic singularities also contribute nontrivially to this correction term, and so for simplicity, we assume throughout that $\pi_1(W)$ is torsion free. In this case, Hirzebruch's formula becomes 
\[ \si(\wt{W}) = \de(E_1)+\dots +\de(E_r) \] 
where $E_1,\dots,E_r$ is a complete set of cusp ends of $W$ given from the thick-thin decomposition and $\si(\wt{W})$ denotes the signature of $\wt{W}$. The definition of the terms $\de(E_j)$ are given as follows. Associated to each cusp end is the $\pi_1(W)$--conjugacy class of a maximal peripheral subgroup $\Ga_j$. The group $\Ga_j$ is conjugate in $\PSL(2;k)$ to a subgroup of the familiar form $\tri(M_j,V_j)$. In turn, for the pair $(M_j,V_j)$, we have an associated Shimuzu $L$--function $L(M_j,V_j,s)$---see \cite{Shimizu63}---defined by
\[ L(M,V,s) = \sum_{\be \in (M_j\smin \set{0})/V_j} \frac{\sign(N_{k/\Q}(\be))}{(N_{k/\Q}(\be))^s} \]
where $N_{k/\Q}$ is the norm map. With this, the invariant $\de(E_j)$ is defined to be
\[ \de(E_j) = \frac{-\vol(M_j)}{\pi^2}L(M_j,V_j,1) \]
where $\vol(M_j)$ is the volume of $\R^2/M$ with respect to the pairing $\tr_{k/\Q}$. Equivalently, 
\[ \vol(M_j) = \abs{\det(\be_i^{(j)})}, \] 
where $\be_1,\be_2$ is a $\Z$--module basis for $M_j$ and $\be_i^{(1)}$ and $\be_i^{(2)}$ denote the image of $\be_i$ under the two real embeddings of $k$ into $\R$.

\begin{thm}[Hirzebruch;\cite{Hirzebruch73}]\label{T:Signature}
If $W$ is a Hilbert modular manifold with exactly one cusp, then
\[ \si(\wt{W}) = \frac{-\vol(M)}{\pi^2}L(M,V,1) \]
for the unique $\pi_1(W)$--conjugacy class $\tri(M,V)$. 
\end{thm}

As we seek an integrality condition, it is convenient to change the pair $M,V$. Associated to the $\Z$--module $M$ is the \emph{dual lattice} $M^*$ defined to be the image of $M$ under the duality pairing provided by $\tr_{k/\Q}$.

\begin{prop}\label{P:Dual}
For a horosphere $\Cal{H}$ stabilized by $\tri(M,V)$ and $\tri(M^*,V)$, $\Cal{H}/\tri(M,V)$ and $\Cal{H}/\tri(M^*,V)$ are diffeomorphic $\Sol$ 3--manifolds.
\end{prop}

\begin{pf} Let $\vp_M,\vp_{M^*}\co V \lra \SL(2;\Z)$  be the holonomy representations for $\tri(M,V)$ and $\tri(M^*,V)$. The pairing $\tr_{k/\Q}$ can be viewed as an element of $\la \in\SL(2;\Z)$ such that $\la M = M^*$. By construction $\vp_{M^*} = \la(\vp_M)\la^{-1}$, and so we have an isomorphism $\rho\co \tri(M,V) \lra \tri(M^*,V)$ given by
\[ \rho(\be,\vp_M(\al)) = (\la \be, \la\vp_M(\al)\la^{-1}). \]
The proof is completed by appealing to the smooth rigidity theorem of Mostow \refT{Rig}.
\end{pf}

Hecke \cite{Hecke83} related the $L$--functions $L(M,V,s)$ and $L(M^*,V,s)$ by the functional equation $H(M,V,s) = (-1)^sH(M^*,V,1-s)$, where
\[ H(M,V,s) = \brac{\Ga\pr{\frac{s+1}{2}}}^2\pi^{-(s+1)}\brac{\vol(M)}^sL(M,V,s) \]
The specialization of this functional equation at $s=1$ produces 
\begin{align*}
\pr{\Ga(1)}^2\pi^{-2}\vol(M)L(M,V,1) &= -\pr{\Ga\pr{\frac{1}{2}}}^2\pi^{-1}L(M^*,V,0) \\
L(M^*,V,0) &= -\frac{\vol(M)}{\pi^2}L(M,V,1),
\end{align*}
and thus from this and \refT{Signature}, we obtain 
\begin{equation}\label{e:eta}
\si(\wt{W}) = L(M^*,V,0).
\end{equation}

It is at this point that we take stock in what has been done. For a 1--cusped Hilbert modular manifold $W$ with cusp cross-section $S$, we have associated to $S$ the invariant $\de(S \times \R^+)$. As both $M$ and $V$ depend on the associated $\Sol$ metric on $S$ afforded by its embedding as a cusp cross-section, the invariant $\de(S\times \R^+)$ depends on the associated $\Sol$ metric on $S$. Our goal is to use the integrality of $\si(\wt{W})$ and (\ref{e:eta}) to produce an obstruction for $S$ to topologically occur in this geometric setting. For this, it remains to show the invariant $\de(S \times \R^+)$ is independent of the $\Sol$ structure on $S$.  

Given a peripheral group $\tri(M,V)$ and stabilized horosphere $\Cal{H}$, the metric on $\Hy_\R^2 \times \Hy_\R^2$ endows $\Cal{H}$ with a $\tri(M,V)$--invariant metric $g_{\Cal{H},M,V}$. Consequently the metric $g_{\Cal{H},M,V}$ descends to quotient $\Cal{H}/\tri(M,V)$ and endows $\Cal{H}/\tri(M,V)$ with a complete $\Sol$ structure that depends on the horosphere $\Cal{H}$ only up to similarity. 

The formula (\ref{e:eta}) was also established in \cite{ADS83} (see also \cite{Muller84}) where $L(M^*,V,0)$ was reinterpreted as the $\eta$--invariant of an adiabatic limit.

\begin{thm}[Atiyah--Donnely--Singer;\cite{ADS83}]\label{T:ADS1}
\[ L(M^*,V,0) = \lim_{\eps \lra 0} \eta(\Cal{H}/\tri(M^*,V),g_{\Cal{H},M^*,V}/\eps). \]
\end{thm}

More generally, given any $\Sol$ structure $g$ on $S$, we can define 
\[ \de(S,g) = \lim_{\eps \lra 0} \eta(S,g/\eps). \] 
The last ingredient for proof of \refT{GB} is the independence of $\de(S,g)$ from $g$, a result established by Cheeger and Gromov \cite{CheegerGromov85} (see \cite{Cheeger87} for a treatment specific to $\Sol$).

\begin{thm}[Cheeger--Gromov;\cite{CheegerGromov85}]\label{T:Eta}
$\de(S,g)$ is a topological invariant of the $\Sol$ 3--manifold $S$.
\end{thm}

We are now in position to state and prove the principal observation needed in the proof of \refT{GB} (compare with \cite{LongReid00}).

\begin{thm}\label{T:GeometricBounding}
If $S$ is diffeomorphic to a cusp cross-section of a 1--cusped Hilbert modular manifold, then $\de(S) \in \Z$.
\end{thm}

\begin{pf}
If $(S,g)$ arises as a cusp cross-section of a 1--cusped Hilbert modular manifold $W$, then there is an isometric embedding $f\co (S,g) \lra W$ 
onto a cusp cross-section of $W$. Let $f_*(\pi_1(S)) =\tri(M,V)$ with associated horosphere $\Cal{H}$ selected such that $\Cal{H}/\tri(M,V)$ is embedded in $W$. By \refP{Dual}, $\Cal{H}/\tri(M^*,V)$ is diffeomorphic to $S$, though equipped with the metric $g_{\Cal{H},M^*,V}$. From the computation above in combination with \refT{ADS1}, $\si(\wt{W}) = \de(S,g_{\Cal{H},M^*,V})$ and by \refT{Eta}, the right hand side depends only on the topological type of $S$. Since $\si(\wt{W})$ is in $\Z$, $\de(S)$ is in $\Z$ as asserted.
\end{pf}

\begin{pf}[Proof of \refT{GB}]
To prove \refT{GB}, by \refT{GeometricBounding}, it suffices to find a $\Sol$ 3--manifold $S$ for which $\de(S) \notin \Z$. For $k=\Q(\sqrt{3})$, the standard Hilbert modular surface $W$ over $k$ has precisely one cusp, since the number of cusps of a standard Hilbert modular surface over $k$ is the ideal class number of $k$. Setting $S$ to be an embedding cusp cross-section of $W$, the proof is completed by appealing to \cite{Hirzebruch73}. Specifically, Hirzebruch showed $\de(S) = -1/3$.
\end{pf}

\begin{rem} It is unknown to the author whether or not there exist 1--cusped Hilbert modular manifolds. In addition, the number fields $\Q(\sqrt{6})$, $\Q(\sqrt{21})$ and $\Q(\sqrt{33})$ also have standard Hilbert modular surfaces with precisely one cusp for which the associated invariant $\de(S) \notin \Z$. In each of these cases, $\de(S)=-2/3$ (see \cite[p. 236]{Hirzebruch73}). 
\end{rem} 

\chapter{Examples and low dimensional considerations}

\section{Prime order holonomy}

For an AB-group $\Ga$ modelled on $\Fr{N}_{2n-1}$ with cyclic order $p$ holonomy, $C_p$, where $p$ is an odd prime, the holonomy is necessarily complex and acts trivially on the center of the Fitting subgroup. Taking the quotient of $\Ga$ by its center gives birth to a $(2n-2)$--dimensional Bieberbach group with $C_p$--holonomy. Hence, there exists a faithful representation of $C_p$ into $\GL(2n-2;\Z)$ afforded to us by the Bieberbach theorems. This can occur only when $p-1\leq 2n-2$, and that such AB-groups exist in dimension $2n-1$ can be shown by explicit construction. The following proposition shows the existence of infinitely many AB-groups modelled on $\Fr{N}_{2n-1}$ (for infinitely many $n$) which are not arithmetically admissible. 

\begin{prop}\label{P:411} If $\Ga_p$ is an AB-groups modelled on $\Fr{N}_{2n-1}$ with holonomy $C_p$, $2(n - 1) = p-1$ and $\Ga_p$ is arithmetically admissible, then $p \equiv 3 \mod 4$.
\end{prop}

\begin{pf} If $\Ga_p$ is arithmetically admissible, by \refC{Holonomy}, there exists a faithful representation ($k$ is an imaginary quadratic number field) 
\[ \rho\co C_p \lra \GL\pr{\frac{p-1}{2};k}. \] 
Let $k_\rho$ denote the field generated by the traces of $\rho(\xi)$ for a generator $\xi$ of $C_p$ and note $k_\rho$ is contained in $k$. The representation $\rho$ is conjugate to one which decomposes into a direct sum of characters $\chi_j\co C_p \lra \C^\times$. Each of these characters $\chi_j$ is of the form $\chi_j(\xi) = \zeta_p^{n_j}$. Therefore 
\[ \tr(\rho(\xi)) = \sum_{j=1}^{\frac{p-1}{2}} \zeta_p^{n_j}. \]
Since $\rho$ is faithful, for some $j$, 
\[ n_j \ne 0 \mod p. \] 
By considering the cyclotomic polynomial $\Phi_p(x)$, we deduce $\tr(\rho(\xi))$ is not in $\Q$ and so $k_\rho$ is a nontrivial extension of $\Q$. On the other hand, from the decomposition above, $k_\rho$ is contained in $\Q(\zeta_p)$. As $[k:\Q]=2$, it must be that $k=k_\rho$. Hence $\Q(\zeta_p)$ contains an imaginary quadratic extension of $\Q$. By quadratic reciprocity, this can happen if and only if $p \equiv 3 \mod 4$.
\end{pf}

It is worth noting for $p \equiv 3 \mod 4$, \refC{Holonomy} can be used to show such AB-groups are arithmetically admissible and the field $k$ is the unique imaginary quadratic number field in $\Q(\zeta_p)$. Moreover, the holonomy generator acts by the matrix
$\res_{\Q(\zeta_p)/k}(\zeta_p)$.  

\begin{rem} When $p >5$, we get an obstruction without appealing to \refT{324}. In this case, we have an injection 
\[ \rho\co C_p \lra \Uni(H;k) < \GL\pr{\frac{p+1}{2};k}. \]
So long as $(p+1)/2 < p-1$, this implies $p \equiv 3 \mod 4$.
\end{rem}

In the quaternionic setting, take $\Ga$ modelled on $\Fr{N}_{7}(\BB{H})$ with $C_5$--holonomy where the action of $C_5$ on the center of $\Fr{N}_7(\BB{H})$ is trivial. If $\Ga$ is arithmetically admissible, we obtain an injection of $C_5$ into $A$, a ramified quaternion algebra over $\Q$. However, this is impossible since such $A$ do not contain elements of order five (see \cite{Amitsur55} or \cite{Vigneras80}). Thus $\Ga$ cannot be arithmetically admissible.  

\begin{prop} There exists infinitely many nonarithmetically admissible AB-groups modelled on $\Fr{N}_7(\BB{H})$.
\end{prop} 

\section{Central products}

In the complex setting, we get arithmetically inadmissible examples in every dimension by a gluing constructing we call the \emph{central product}. The details are as follows. For a pair of AB-groups $\Ga_1,\Ga_2$ modelled on $\Fr{N}_{2n_1-1}$ and $\Fr{N}_{2n_2-1}$, respectively, we define 
\[ \Ga_1 \times_{c} \Ga_2 = \Ga_1 \times \Ga_2 / N_{\Ga_1\times \Ga_2}(\innp{c_1c_2^{-1}}), \] 
where $c_j$ is the generator of the center of the Fitting subgroup $L_j$ of $\Ga_j$ and $N_{\Ga_1\times \Ga_2}(c_1c_2^{-1})$ denotes the normal closure of $\innp{c_1c_2^{-1}}$ in $\Ga_1\times\Ga_2$. We call this the \emph{central product} of $\Ga_1$ and $\Ga_2$. The following lemma shows under some mild assumptions, the central product is also an AB-group modelled on $\Fr{N}_{2(n_1+n_2)-1}$.

\begin{lemma}
If the holonomy groups $\te_1,\te_2$ are either both complex or both anticomplex, then $\Ga_1\times_c \Ga_2$ is an AB-group modelled on $\Fr{N}_{2(n_1+n_2)-1}$.
\end{lemma}

\begin{pf} As we only wield this when the holonomy groups are complex, the anticomplex case is left for the reader. To begin, we have natural inclusions of $\Fr{N}_{2n_1-1}$ and $\Fr{N}_{2n_2-1}$ into $\Fr{N}_{2(n_1+n_2)-1}$ induced by
\[ \rho_j\co \C^{n_j} \lra \C^{n_1} \op \C^{n_2}. \]
In fact, this yields inclusions of 
\[ \C^{n_j} \rtimes \GL(n_j;\C) \lra \C^{n_1+n_2} \rtimes \GL(n_1+n_2;\C). \] 
As a result, we have an injective homomorphism of $\Ga_1 \times_c \Ga_2$ onto an AB-group in $\Aff(\Fr{N}_{2(n_1+n_2)-1})$. By selecting the maps $\rho_j$, we can ensure the induced maps agree on the center of the Fitting subgroups of $\Ga_1$ and $\Ga_2$. With this selection, the induced map becomes an isomorphism of $\Ga_1 \times_c \Ga_2$ with an AB-group in $\Fr{N}_{2(n_1+n_2)-1}$, as desired.
\end{pf} 

Using central products we can construct many arithmetically inadmissible AB-groups. We summarize this in the following theorem which proves \refT{Ob} (a).

\begin{thm}[Central product theorem]\label{T:ConGlue}
If $\Ga_1$ and $\Ga_2$ be AB-groups modelled on $\Fr{N}_{2n_1-1}$ and $\Fr{N}_{2n_2-1}$, defined over $k_1$ and $k_2$, respectively, and  $\Ga_1$ and $\Ga_2$ are both complex or anticomplex, then
\begin{itemize}
\item[(a)] $\Ga_1 \times_c \Ga_2$ is $k_1k_2$--defined,
\item[(b)] $\Ga_1 \times_c \Ga_2$ is $k$--arithmetic if and only if both $\Ga_1$ and $\Ga_2$ are $k$--arithmetic, and
\item[(c)] there exist arithmetically inadmissible AB-groups modelled on $\Fr{N}_{2n-1}$ for all $n \geq 3$.
\end{itemize}  
\end{thm}

\begin{pf} For (a), let $M_1(X)$ and $M_2(X)$ denote the maximal compact groups defined over $k_1$ and $k_2$ for which $\Ga_j$ injects into the $k_j$--points of $\Fr{N}_{2n_j-1}(X) \rtimes M_j(X)$. Associated to each of these maximal compact groups is a $k_j$--defined hermitian form $H_j$. Let $H = H_1 \op H_2$ and $M(X)$ be the associated maximal compact subgroup. $M(X)$ is $k_1k_2$--define and we have an injection (into the $k_1k_2$--points) 
\[ \rho\co \Ga_1\times_c \Ga_2 \lra \Fr{N}_{2(n_1+n_2)-1}(X) \rtimes M(X). \] 
Thus $\Ga_1\times_c\Ga_2$ is $k_1k_2$--defined.  

For (b), if $\Ga=\Ga_1 \times_c \Ga_2$ is $k$--arithmetic, then $\Ga$ is isomorphic to a maximal peripheral subgroup of an arithmetic lattice $\La$ in $\Uni(H)$, where $H$ is a signature $n_1+n_2-1$ hermitian form defined over an imaginary quadratic number field
$k$. It must be that $\Ga_j$ injects into a subgroup $\La_j$ of $\La$ which is maximal with respect to stabilizing a $k$--defined
complex subspace $\C^{n_j,1}$. As $\La_j$ is an arithmetic lattice in a smaller isometry group (whose model form is the restriction of $H$ to
the complex subspace $\C^{n_j,1}$), this implies $\Ga_j$ is $k$--arithmetic for $j=1,2$. The reverse implication follows immediately from (a).

For (c), by \refP{411}, there exists an arithmetically inadmissible AB-group modelled on $\Fr{N}_5$. To obtain examples in higher
dimensions, we take central products of this example with other AB-groups and apply (b).
\end{pf} 

The quaternionic setting can be handled similarly. We construct examples in every dimension by taking central products with the inadmissible example $\Ga$ in $\Fr{N}_7(\BB{H})$ given above.

\begin{cor}\label{C:64} There exist arithmetically inadmissible AB-groups modelled on $\Fr{N}_{4n-1}(\BB{H})$ for all $n \geq 2$.
\end{cor}

\section{$\Nil$ 3--manifolds and the proof of \refT{Nil3Thm}}

Fix a closed $\Nil$ 3--manifold $M$ with $\pi_1(M)=\Ga$. By the generalized Bieberbach theorems (see \cite{Dekimpe96}), to prove \refT{Nil3Thm}, it suffices to show $\Ga \cong \tri(\La)$, where $\La$ is an arithmetic lattice in $\Isom(\Hy_\C^2)$. In fact, by \refT{Cusp} it suffices to construct an injective homomorphism $\vp\co \Ga \lra \tri(\La)$. To this end, let $\Fr{N}(3) = \Fr{N}_3\rtimes \Uni(1)$ and $\iota$ be the isometry of $\Hy_\C^2$ induced by conjugation. For a subring $R \su \C$, we define $\Fr{N}_3(R) = R \times \Ima R$ with the induced group operation and set $\Fr{N}(3,R) = \Fr{N}_3(R) \rtimes \Uni(1;R)$. For the statement of our next result, let $\zeta_3$ be a primitive third root of unity, say $\zeta_3=-1/2+\sqrt{-3}/2$.  

\begin{thm}\label{T:22} If $M$ is a closed $\Nil$ 3--manifold $M$ and $\Ga =\pi_1(M)$, then there exists a faithful representation 
\[ \vp\co \Ga \lra \innp{\Fr{N}(3,\Cal{O}_k),\iota} \] 
with $k=\Q(i)$ or $\Q(\zeta_3)$.
\end{thm}

\begin{pf} We begin by summarizing the strategy of the proof, which depends heavily on the list of presentations for the fundamental group of a closed $\Nil$ 3--manifold found in the appendix. The idea is to show an injective homomorphism on the Fitting subgroup $\innp{a,b,c}$ can be promoted to the full 3--manifold group (see \cite[Thm 3.1.3]{Dekimpe96}). To get a representation with the coefficients in $\Z[i]$ or $\Z[\zeta_3]$, we are reduced to solving some simple equations. The details are as follows. 

In the lemma below, let $p_1\co \Fr{N}_3 \lra \C$ be projection onto the first factor. 

\begin{lemma}\label{L:311}
\begin{itemize}
\item[(a)] If $a$ and $b$ are as above, 
\[ \rho\co \innp{a,b,c} \lra \Fr{N}_3 \]
is a homomorphism, $p_1(\rho(a))$ and $p_1(\rho(b))$ are $\Z$--linearly independent, and $c\notin \ker \rho$, then $\rho$ is injective. 
\item[(b)] In addition, if 
\[ \rho^{-1}(\rho(\innp{a,b,c})) = \innp{a,b,c}, \] 
and $\rho_{|\innp{a,b,c}}$ is an injective homomorphism, then $\rho$ is an injective homomorphism. 
\end{itemize}
\end{lemma}

\begin{pf} For (a), let $w$ be in $\ker \rho$, write
\[ w = a^{n_1}b^{n_2}c^{n_3}, \] 
and set 
\[ \rho(a) = (v_1,t_a), \quad \rho(b)=(v_2,t_b). \] 
Since $[a,b]=c^k$, it must be that $\rho(c)=(0,s)$ as $[\Fr{N}_3,\Fr{N}_3]=\set{(0,t)~:~t \in \R}$. With this we have 
\[ \rho(w) = (n_1v_1+n_2v_2,n_1t_a+n_2t_b+2\Ima\innp{n_1v_1,n_2v_2}+n_3s). \] 
Since $w$ is in $\ker \rho$, 
\[ n_1v_1 + n_2v_2 = 0. \] 
The $\Z$--linear independence of $v_1$ and $v_2$ implies $n_1=n_2=0$. Therefore $n_3s=0$, and so $n_3=0$, as $s \ne 0$.

For (b), let $w$ be in $\ker \rho$ and write
\begin{equation}\label{word}
w = a^{n_1}b^{n_2}c^{n_3}\al^{s_1}\be^{s_2}, \quad n_1,n_2,n_3,s_1,s_2\in \Z.
\end{equation}
Using this form of $w$, we have
\[ \rho(a^{n_1}b^{n_2}c^{n_3}) = \rho(\al^{-s_1}\be^{-s_2}). \]
By assumption 
\[ \rho^{-1}(\rho(\innp{a,b,c}))=\innp{a,b,c}, \] 
and so $\al^{s_1}\be^{s_2}$ is contained in $\innp{a,b,c}$. This in tandem with (\ref{word}) implies $w$ is in $\innp{a,b,c}$. However $\rho_{|\innp{a,b,c}}$ is one-to-one, and so $w=1$.
\end{pf}

For (b), we only require that if $\rho(\al^{s_1}\be^{s_2}) \in \rho(\innp{a,b,c})$ then $\al^{s_1}\be^{s_2} \in \innp{a,b,c}$. Indeed, we need only check for $s_1 \in \set{1,\dots,k_\al}$ and $s_2 \in \set{1,\dots,k_\be}$ where $k_\eps$ is the first integer such that $\eps^{k_\eps}\in\innp{a,b,c}$, $\eps = \al$ or $\be$.

Let $L=\innp{a,b,c}$ and define two homomorphisms
\[ \vp_3,\vp_4\co L \lra \Fr{N}_3 \] by 
\[ \vp_j(a) = (1,0), \quad \vp_j(b) = (\zeta_j,0). \] 
This determines $c$, since some power of $c$ is a commutator of $a$ and $b$. By \refL{311} (a), both maps are injective homomorphisms. 

We extend this to $\Ga$ by declaring
\begin{equation}\label{E:311}
\vp_j(\al) = (z_1,t_1,\eta_1), \quad \vp_j(\be) = (z_2,t_2,\eta_2),
\end{equation}
where $z_1,z_2\in\C$, $t_1,t_2\in\R$, and $\eta_1,\eta_2\in\innp{\Uni(1),\iota}$. 

To solve \refE{311}, we simply use the presentations in the appendix to ensure this yields a homomorphism. By applying \refL{311} (b), one can see these solutions yield injective homomorphisms. For clarity, we solve the equations for the second family (2) and give a list of the equations and solutions for the seventh family (7).

The second family has presentation
\begin{align*}
<a,b,c,\al~:~ & [b,a]=c^k,~[c,a]=[c,b]=[c,\al]=1,\\
&~\al a = a^{-1}\al, \al b = b^{-1} \al, \al^2 = c>
\end{align*}
with $k \in 2\N$. For this family we take the map $\vp_4$. First, consider the relation $[b,a]=c^k$. Under $\vp_4$, the left hand side becomes
\begin{align*}
[\vp_4(b),\vp_4(a)] &= [(i,0,1),(1,0,1)] \\
&= (i+1,2\Ima\innp{i,1},1)(-i-1,2\Ima\innp{-i,-1},1) \\
&= (0,4,1).
\end{align*}
Since $[\vp_4(b),\vp_4(b)]=\vp_4(c)^k$, it follows that $\vp_4(c)=(0,4/k,1)$.

Next, consider the relation $\al^2 = c$. Under $\vp_4$, we have
\begin{align*}
\vp_4(\al^2) &= (z_1,t_1,\eta_1)(z_1,t_1,\eta_1) \\
&= (z_1 + \eta_1z_1,2t_1 + 2\Ima\innp{z_1,\eta_1z_1},\eta_1^2) \\
&= (0,4/k,1) \\
&=\vp_4(c).
\end{align*}
In particular, $\eta_1^2=1$. If $\eta_1=1$, the above injection would yield an isomorphism between a group in the first family with a group in the second family. This is impossible, therefore $\eta_1=-1$. By considering the first coordinate equation with $\eta_1=-1$, we get no  information. The second coordinate equation is $2t_1 = 4/k$, therefore $t_1=2/k$. One can now check that $[c,\al]=1$, regardless of $z_1$. 

Now, we take the relation $\al a = a^{-1} \al$. We have
\begin{align*}
\vp_4(\al a) &= (z_1,2/k,-1)(1,0,1) = (z_1 - 1, 2/k - 2\Ima z_1, -1).
\end{align*}
On the other hand,
\begin{align*}
\vp_4(a^{-1}\al) &= (-1,0,1)(z_1,2/k,-1) = (-1+z_1,2/k + 2\Ima z_1,-1).
\end{align*}
The first and last coordinate equations yield no information, while the second coordinate equation yields $4\Ima z_1=0$. Hence $\Ima
z_1=0$. 

Lastly, we have the relation $\al b = b^{-1}\al$. We have
\begin{align*}
\vp_4(\al b) &= (z_1,2/k,-1)(i,0,1) = (z_1 - i,2/k + 2\RE z_1,-1).
\end{align*}
On the other hand,
\begin{align*}
\vp_4(b^{-1}\al) &= (-i,0,1)(z_1,2/k,-1) = (-i+z_1,2/k - 2\RE z_1,-1).
\end{align*}
As above, the first and second coordinates yields no information, while the second coordinate implies $\RE z_1=0$.

Hence, we deduce from the above computations, the desired homomorphism $\rho$ is defined by
\begin{align*}
\rho(a) &= (1,0,1), & \rho(b) = (i,0,1) & \\
\rho(c) &= (0,4/k,1) & \rho(\al) = (0,2/k,-1). & 
\end{align*}

Visibly, these solutions are in $\Q(i)$ and not $\Z[i]$. This is rectified by conjugating the above representation by a dilation of $2k$. This dilation is linear on the first factor, quadratic on the second factor, and trivial on the third factor. The resulting faithful representation is
\begin{align*}
\rho(a) &= (2k,0,1), & \rho(b) = (2ki,0,1) & \\
\rho(c) &= (0,16k,1) & \rho(\al) = (0,8k,-1). & 
\end{align*}
The faithfulness of $\rho$ follows from \refL{311} (b) upon verifying the conditions of this lemma are met. The injectivity of $\rho_{|\innp{a,b,c}}=\vp_4$ follows from \refL{311} (a). To check 
\[ \rho(\rho^{-1}(\innp{a,b,c}))=\innp{a,b,c}, \] 
by the remark proceeding the proof of \refL{311}, it suffices to show $\al$ is not in $\rho(\rho^{-1}(\innp{a,b,c}))$, the validity of which obvious.  
 
For the seventh family (7), we have the presentation
\begin{align*}
<a,b,c,\al~:~&[b,a]=c^k,~[c,a]=[c,b]=[c,\al]=1,~\al a =
ab\al,\\
&~\al b=a^{-1}\al,~\al^6=c^{k_1}>
\end{align*}
with
\begin{align*}
k \equiv 0 \mod 6,~k_1 =1, & \text{ or } k \equiv 4 \mod 6,~k_1=1,
\text{ or }\\
k \equiv 0 \mod 6,~k_1 =5, & \text{ or } k \equiv 2 \mod 6,~k_1=5.
\end{align*}

We take $\vp_3$ in this case. By considering all the relations, we get
\begin{align*}
(0,4/k,1) &= (0,s,1)\\
(z_1+\eta_1,t_1+2\Ima\innp{z,\eta_1},\eta_1) &= (1+\zeta_3 + z_1, t_1+  2\Ima\innp{1+\zeta_3,z_1} \\
& + 2\Ima\innp{1,\zeta_3} + t_1,\eta_1) \\
(z_1 + \eta_1\zeta_3,t_1+2\Ima\innp{z_1,\eta_1\zeta_3},\eta_1) &=
(-1+z_1,t_1+2\Ima\innp{-1,z_1},\eta_1).
\end{align*}
We omit the last relation $\al^6 =c^{k_1}$ as it is quite long. Note the commutator relations (aside from $[b,a]=c^k$) are all trivially satisfied. Solving these equations and conjugating by a dilation of $12k$ to get the coefficients in $\Z[\zeta_3]$, we have
\begin{align*}
\vp_3(a) &= (12k,0,1), & \vp_3(b)= (12k\zeta_3,0,1), \\ \vp_3(c) &=
(0,288k\sqrt{3},1), & \vp_3(\al)= (-6k,12k(4k_1+3k)\sqrt{3},\zeta_6). &
\end{align*}
That this is faithful again follows from \refL{311} and the remark following its proof.  
 
The other remaining families are handled similarly and in the appendix, we provide a list of the matrices obtained from this venture and the post-composition of the resulting representation with $\psi$ (see below).
\end{pf}

The proof of \refT{Nil3Thm} is completed by following the injection $\rho$ from \refT{22} with the injection
\[ \psi\co \Fr{N}_3 \rtimes \innp{\Uni(1),\iota} \lra \Uni(2,1) \]
given by
\[ (\xi,t,\iota^\eps U) \lmto \begin{pmatrix} 1 & \xi & \xi \\ -\ol{\xi} &
1-\frac{1}{2}(\norm{\xi}^2 -it) & -\frac{1}{2}(\norm{\xi}^2-it) \\
\ol{\xi} & \frac{1}{2}(\norm{\xi}^2-it) & 1+\frac{1}{2}(\norm{\xi}^2-it)
\end{pmatrix} \begin{pmatrix} U & 0 & 0 \\ 0 & 1 & 0 \\ 0 & 0 & 1
\end{pmatrix} \iota^\eps. \] 

\section{Infranil 5--manifolds}

Using the list of isomorphism types of holonomy groups for infranil 5--manifolds given in \cite{DekimpeEick02}, we can carry out the same analysis. 

\begin{prop} The only complex holonomy groups which yield arithmetically inadmissible groups are $C_5$, $C_{10}$, $C_{12}$ and $C_{24}$.
\end{prop}

To see this, by \refC{Holonomy}, it suffices to check the field of definition for all the distinct representations of the holonomy group in $\GL(2;\C)$. For most of the groups, every representation will be conjugate to one defined over an imaginary quadratic number field. Both $C_{12}$ and $C_{24}$ can arise via central products of $\Nil$ 3--manifold groups for which \refT{ConGlue} can be applied to show the resulting AB-groups are arithmetically inadmissible. For $C_5$--holonomy, we can apply \refP{411}, which in turn yields the result for $C_{10}$--holonomy, since any arithmetic representation for an AB-group with $C_{10}$--holonomy would yield one for an AB-group with $C_5$--holonomy.

\section{$\Sol$ 3--manifolds}

In \cite{Scott83}, Scott proved every $(2,1)$--torus bundles admits either a Euclidean, $\Nil$, or $\Sol$ structure. The following essentially establishes \refT{Sol3Thm}. 

\begin{prop}\label{P:Scott}
If $M$ is an orientable $(2,1)$--torus bundle which admits a $\Sol$ structure, then there exists a faithful representation 
\[ \rho\colon \pi_1(M) \lra \Cal{O}_k \rtimes \Cal{O}_k^\times \] 
for some real quadratic number field $k$.
\end{prop}

\begin{pf}
For any $(2,1)$--torus bundle $M$, let the $\Z$--action be given by 
\[ A = \begin{pmatrix} a & b \\ c & d \end{pmatrix}. \]
If the order of $A$ is finite, then $\pi_1(M)$ is a Bieberbach group and $M$ admits a Euclidean structure. Therefore we assume the
$A$ is not a torsion element. If $A$ is not diagonalizable, then some power of $A$ is conjugate to 
\[ \begin{pmatrix} 1 & \al \\ 0 & 1 \end{pmatrix} \] 
with $\al \ne 0$. In this case, $M$ admits a $\Nil$ structure. Thus, we assume $A$ is diagonalizable, and in this case we have \[ \begin{pmatrix} \be & 0 \\ 0 & \be^{-1} \end{pmatrix} \] 
for a conjugate of $A$. It follows, since $A \in \SL(2;\Z)$, that $\be$ and $\be^{-1}$ are algebraic integers in the real quadratic field $\Q(\be)$. Thus the representation 
\[ \vp\colon \Z \lra \GL(2;\Z) \]
is conjugate to $\res_{k/\Q}(\chi)$, where $\chi\colon \Z \lra \Cal{O}_k^\times$ is given by $\chi(1) = \be$. Therefore by the remark following \refT{Class}, we have a faithful representation 
\[ \rho\colon \pi_1(M) \lra \Cal{O}_k \rtimes \Cal{O}_k^\times, \] 
as asserted.
\end{pf}

Via \refP{Scott}, every $\Sol$ 3--manifold group does faithfully represent into $\Isom((\Hy_\R^2)^2)$. Those arising as cusp cross-sections of Hilbert modular surfaces are precisely the ones whose fundamental group faithfully represents into the identity component 
of $\Isom((\Hy_\R^2)^2)$. However, the quotients of those groups which fail to map into the identity component do produce finite volume quotients which possess 2--fold covers that are Hilbert modular surfaces. For this reason, we call such quotients \emph{generalized Hilbert modular varieties}. Given this, \refT{Sol3Thm} follows from this discussion in combination with \refT{Torus}.

\section{Two $\Sol$ examples}

The following example illustrates the ideas in the proofs of \refT{Class} and \refT{Sol3Thm}. 

\begin{exa} Let $M$ be a $(2,1)$--torus bundle with $\Z$--action given by 
\[ A = \begin{pmatrix} 1 & 2 \\ 1 & 3 \end{pmatrix}. \]
With this action, $\pi_1(M)$ has a presentation of the form (we are assuming $\pi_1(M)$ is a split extension which is always the case; see \cite{Scott83})
\[ \innp{a_1,a_2,b~:~ [a_1,a_2]=1,~ba_1 = a_1a_2b,~ba_2=a_1^2a_2^3b}. \]
To obtain a faithful representation 
\[ \rho\colon \pi_1(M) \lra \Cal{O}_k \rtimes \Cal{O}_k^\times \] 
for some quadratic number field $k$, we first compute the eigenvalues of $A$. The characteristic polynomial for $A$ is 
\[ c_A(t)=t^2-4t+1, \] 
which has roots $2 \pm \sqrt{3}$. Let $k=\Q(\sqrt{3})$ and write  
\[ \rho(a_1) = (x_1+y_1\sqrt{3},1), \quad \rho(a_2) = (x_2+y_2\sqrt{3},1), \quad \rho(b)= (0,2+\sqrt{3}) \]
where $x_1,x_2,x_3,$ and $x_4$ are to be determined. Using the presentation above, we are now reduced to solving a system of equations in $x_1,x_2,x_3$, and $x_4$ to ensure $\rho$ is an injective homomorphism. By construction $\rho([a_1,a_2]) = 1$. The other two relations yield the equations
\begin{align*}
(2+\sqrt{3})(x_1+y_1\sqrt{3}) &= x_1+y_1\sqrt{3}+x_2+y_2\sqrt{3} \\
(2+\sqrt{3})(x_2+y_2\sqrt{3}) &= 2(x_1+y_1\sqrt{3})+3(x_2+y_2\sqrt{3}).
\end{align*} 
Solving, we get the faithful representation
\[ \rho(a_1) = (1,1), \quad \rho(a_2) = (1+\sqrt{3},1), \quad \rho(b) = (0,2+\sqrt{3}). \]
\end{exa}

Let $k$ be a totally real, cubic Galois extension of $\Q$ and $\be \in \Cal{O}_{k,+}^\times$ be of infinite order. By restricting scalars, we can view $\be$ as an element of $\SL(3;\Z)$. Let $M$ be a $(3,1)$--torus bundle with $\Z$--action given by $\be$. Following the proof of \refT{Class}, we can construct a faithful representation of $\pi_1(M)$ into $\Cal{O}_k \rtimes \Cal{O}_{k,+}^\times$. On
the other hand, any Hilbert modular group defined over $k$ cannot contain $\pi_1(M)$ as a finite index subgroup, since $\rank
\Cal{O}_k^\times = 2$.

The following example is a specific case of the above.

\begin{exa} Define 
\[ A = \begin{pmatrix} 0 & 0 & 1 \\ 1 & 0 & 2 \\ 0 & 1 & -1 \end{pmatrix} \] 
and let $M$ be the $(3,1)$--torus bundle with $\Z$--action given by $A$. The characteristic polynomial for $A$ is $t^3+t^2-2t-1$. This polynomial is irreducible over $\Q$ and has totally real cubic splitting field. Thus, $M$ is a $(3,1)$--torus bundle which is not diffeomorphic to a cusp cross-section of any generalized Hilbert modular variety. However, for the splitting field $k$ of $A$, we do have an injection 
\[ \rho\colon \pi_1(M) \lra \tri(\PSL(2;\Cal{O}_k)). \]
\end{exa} 

\chapter{Presentation of $\Nil$ 3--manifold groups}

The following is a complete list of closed $\Nil$ 3--manifold groups (see \cite[p. 159--166]{Dekimpe96}).

\begin{itemize}
\item[(1)]
\[ <a,b,c~:~[b,a] = c^k,~[c,a] = [c,b] = 1>, \text{ with }k \in \N. \]
\item[(2)]
\begin{align*}
<a,b,c,\al~:~&[b,a]= c^k,~[c,a] = [c,b] = [\al,c] = 1,~\al a = a^{-1} \al \\ & ~\al b = b^{-1}\al,~\al^2 = c>, \text{ with }k \in 2\N.
\end{align*}
\item[(3)]
\begin{align*}
<a,b,c,\al~:~&[b,a] = c^{2k},~ [c,a] = [c,b] = [a,\al] =  1, ~\al c = c^{-1}\al,\\ & ~\al b = b^{-1}\al c^{-k},~\al^2 = a>, \text{ with }k \in \N.
\end{align*}
\item[(4)]
\begin{align*}
<a,b,c,\al,\be~:~& [b,a]=c^{2k},~ [c,a]=[c,b]=[c,\al]=[a,\be]=1,\\ &~\be c= c^{-1}\be,~ \al a = a^{-1} \al c^k, \al b = b^{-1}\al c^{-k},\\ &~ \al^2 = c,~ \be^2 = a, ~ \be b = b^{-1}\be c^{-k}, \\
&~ \al\be = a^{-1}b^{-1}\be \al c^{-k-1}>, \text{ with }k \in \N. 
\end{align*} 
\item[(5)]
\begin{align*}
<a,b,c,\al~:~&[b,a] = c^k,~[c,a]  = [c,b] = [c,\al] = 1,~ \al a  = b\al,\\ &~\al b  = a^{-1}\al,~\al^4  = c^p>,
\end{align*}
$k \in 2\N$ and $p=1$ or $k\in 4\N$ and $p=3$.
\item[(6)]
\begin{align*}
<a,b,c,\al~:~&[b,a] = c^k,~ [c,a] = [c,b] = [c,\al] = 1,~ \al a = b\al c^{k_1},\\ &~\al b = a^{-1}b^{-1}\al,~ \al^3 = c^{k_2}>
\end{align*}
with $k>0$ and
\begin{align*}
k &\equiv 0 \mod 3, \quad k_1=0,~k_2 = 1, \text{ or } \\
k &\equiv 0 \mod 3, \quad k_1=0,~k_2 = 2, \text{ or } \\
k &\equiv 1,2 \mod 3, \quad k_1=1,~k_2=1.
\end{align*}
\item[(7)]
\begin{align*}
<a,b,c,\al~:~&[b,a] = c^k,~ [c,a] = [c,b] = [c,\al] = 1,~\al a = ab\al,\\ &~\al b = a^{-1}\al,~\al^6 = c^{k_1}>,
\end{align*}
with $k>0$ and
\[ k \equiv 0 \mod 6, \quad k_1 = 1,
\text{ or }
k \equiv 4 \mod 6, \quad k_1 = 1, \]
or
\[ k \equiv 0 \mod 6, \quad k_1 = 5,
\text{ or } 
k \equiv 2 \mod 6, \quad k_1 = 5. \]
\end{itemize} 

\chapter{Solutions to \refE{311}}

Below are the resulting matrices obtained by solving \refE{311} and post-composing with the representation $\psi$.

(1)
\begin{align*}
a &= \begin{pmatrix} 1 & 2k & 2k \\ -2k & 1- 2k^2 &
-2k^2 \\ 2k & 2k^2 & 1 + 2k^2 \end{pmatrix}, & b = \begin{pmatrix} 
1 & 2ki & 2ki \\ 2ki & 1-2k^2 & 2k^2 \\ -2ki & 2k^2 & 1+2k^2 \end{pmatrix} & \\
c &= \begin{pmatrix} 1 & 0 & 0 \\ 0 & 1 +8ki & 8ki \\ 0 & -8ki & 1-8ki
\end{pmatrix}.
\end{align*}

(2)
\begin{align*}
a &= \begin{pmatrix} 1 & 2k & 2k \\ -2k & 1- 2k^2 &
-2k^2 \\ 2k & 2k^2 & 1 + 2k^2 \end{pmatrix}, &
b = \begin{pmatrix} 1 & 2ki & 2ki \\ 2ki & 1-2k^2 &
2k^2 \\ -2ki & 2k^2 & 1+2k^2 \end{pmatrix}, & \\
c &= \begin{pmatrix} 1 & 0 & 0 \\ 0 & 1 +8ki & 8ki \\ 0 & -8ki & 1-8ki
\end{pmatrix}, &
\al = \begin{pmatrix} -1 & 0 & 0 \\ 0 & 1+ki & ki \\ 0 & -ki & 1-ki
\end{pmatrix}. &
\end{align*} 

(3)
\begin{align*}
a &= \begin{pmatrix} 1 & 4k & 4k \\ -4k & 1- 8k^2 &
-8k^2 \\ 4k & 8k^2 & 1 + 8k^2 \end{pmatrix}, &
b = \begin{pmatrix} 1 & 4ki & 4ki \\ 4ki & 1-8k^2 &
8k^2 \\ -4ki & 8k^2 & 1+8k^2 \end{pmatrix}, & \\
c &= \begin{pmatrix} 1 & 0 & 0 \\ 0 & 1 +16ki & 16ki \\ 0 & -16ki &
1-16ki \end{pmatrix}, &
\al = \begin{pmatrix} 1 & 2k & 2k \\ 2k & 1 - 2k^2 & -2k^2 \\ -2k &
2k^2 & 1 + 2k^2 \end{pmatrix}\iota. &
\end{align*}

(4)
\begin{align*}
a &= \begin{pmatrix} 1 & 4k & 4k \\ -4k & 1- 8k^2 &
-8k^2 \\ 4k & 8k^2 & 1 + 8k^2 \end{pmatrix},\\
b &= \begin{pmatrix} 1 & 4ki & 4ki \\ 4ki & 1-8k^2 &
8k^2 \\ -4ki & 8k^2 & 1+8k^2 \end{pmatrix}, \\
c &= \begin{pmatrix} 1 & 0 & 0 \\ 0 & 1 +16ki & 16ki \\ 0 & -16ki & 1-16ki
\end{pmatrix}, \\
\al &= \begin{pmatrix} -1 & 2k + 2ki & 2k+2ki \\ -2k+2ki & 1-4k^2 &
-4k^2 \\ 2k-2ki & 4k^2 & 1+4k^2 
\end{pmatrix},  \\
\be &= \begin{pmatrix} 1 & 2k & 2k \\ -2k & 1-4k^2 & -4k^2 \\ 2k &
4k^2 & 1+4k^2 \end{pmatrix} \iota.
\end{align*}

(5)
\begin{align*}
a &= \begin{pmatrix} 1 & 2k & 2k \\ -2k & 1- 2k^2 &
-2k^2 \\ 2k & 2k^2 & 1 + 2k^2 \end{pmatrix}, &
b = \begin{pmatrix} 1 & 2ki & 2ki \\ 2ki & 1-2k^2 &
2k^2 \\ -2ki & 2k^2 & 1+2k^2 \end{pmatrix}, & \\
c &= \begin{pmatrix} 1 & 0 & 0 \\ 0 & 1 +8ki & 8ki \\ 0 & -8ki & 1-8ki
\end{pmatrix}, &
\al = \begin{pmatrix} i & 0 & 0 \\ 0 & 1+\frac{pk}{2}i &
\frac{pk}{2}i \\ 0 & -\frac{pk}{2}i & 1-\frac{pk}{2}i
\end{pmatrix}. &
\end{align*}

(6)
\begin{align*}
a &= \begin{pmatrix} 1 & 24k & 24k \\ -24k & 1 - 288k^2 & -288k^2
\\ 24k & 288k^2 & 1+ 288k^2 \end{pmatrix}, \\
b &= \begin{pmatrix} 1 & -12 + 12\sqrt{3}i & -12 + 12\sqrt{3}i \\
12 + 12\sqrt{3}i & 1 -288k^2 & -288k^2 \\ -12-12\sqrt{3}i & 288k^2 &
1+288k^2 \end{pmatrix},  \\
c &= \begin{pmatrix} 1 & 0 & 0 \\ 0 & 1 + 144k\sqrt{3}i &
144k\sqrt{3}i \\ 0 & -144k\sqrt{3}i & 1-144k\sqrt{3}i \end{pmatrix}, \\
\al &= \begin{pmatrix} 1 & \mu &
\mu \\ 6[k+2k_1]+6\sqrt{3}[k-2k_1]i &
1-\si & -\si \\ -6[k+2k_1]-6\sqrt{3}[k-2k_1]i & \si & 1 +\si
\end{pmatrix} \begin{pmatrix} \zeta_3 & 0 & 0 \\ 0 & 1 & 0 \\ 0 & 0 &
1 \end{pmatrix},
\end{align*}
where
\[ \si = \frac{1}{2}\brac{36(k+2k_1)^2 + 108(k-2k_1)^2 -
 192k\sqrt{3}(k\norm{z}^2 + 2k_2)i} \]
and $\mu = -6[k+2k_1] + 6\sqrt{3}[k-2k_1]i$.

(7)
\begin{align*}
a &= \begin{pmatrix} 1 & 12k & 12k \\ -12k & 1-144k^2 & -144k^2 \\ 12k
& 144k^2 & 1+144k^2 \end{pmatrix}, \\
b &= \begin{pmatrix} 1 & -6k + 6k\sqrt{3}i & -6k + 6k\sqrt{3}i \\ 6k +
6k\sqrt{3}i & 1-144k^2 & -144k^2 \\ -6k - 6k\sqrt{3}i & 144k^2 &
1+144k^2 \end{pmatrix}, \\
c &= \begin{pmatrix} 1 & 0 & 0 \\ 0 & 1 + 288k\sqrt{3}i &
288k\sqrt{3}i \\ 0 & -288k\sqrt{3}i & 1-288\sqrt{3}i \end{pmatrix}, \\
\al &= \begin{pmatrix} 1 & - 6k & -6k \\ 3k-3k\sqrt{3}i & 1 -\chi &
-\chi\\ -3k+3k\sqrt{3}i & \chi & 1+ \chi \end{pmatrix}\begin{pmatrix}
\zeta_6 & 0 & 0 \\ 0 & 1 & 0 \\ 0 & 0 & 1 \end{pmatrix}, 
\end{align*}
where $\chi =  -36k^2 + 12k(4k_1+3k)\sqrt{3}i$.

\bibliographystyle{plain}

\def\cprime{$'$} \def\lfhook#1{\setbox0=\hbox{#1}{\ooalign{\hidewidth
  \lower1.5ex\hbox{'}\hidewidth\crcr\unhbox0}}} \def\cprime{$'$}
  \def\cprime{$'$}

\end{document}